\documentclass[reqno,10pt]{amsart}
\usepackage{latexsym}
\usepackage{amssymb}
\usepackage{amsmath,amsfonts,amsthm,amssymb}
\usepackage{indentfirst}
\usepackage{paralist}
\usepackage{graphics} 
\usepackage{epsfig} 
\usepackage[colorlinks=true]{hyperref}
\usepackage[english]{babel}
\usepackage{amsmath,amssymb}
\hypersetup{urlcolor=red, citecolor=blue}
\newtheoremstyle{theorem}
  {10pt}          
  {10pt}  
  {\sl}  
 {}
  {\bf}  
  {. }    
  { }    
  {}     
\theoremstyle{theorem}

\newtheorem{theorem}{Theorem}[section]
\newtheorem{corollary}{Corollary}[section]

 \newtheorem{lemma}{Lemma}[section]
 
 \newtheorem{remark}{Remark}[section]
  
\numberwithin{equation}{section}

\newtheoremstyle{defi}
  {10pt}          
  {10pt}  
  {\rm}  
  {}  
  {\bf}  
  {. }    
  { }    
  {}     
\theoremstyle{defi}

\allowdisplaybreaks

\begin{document}
\baselineskip = 20 pt
\title[Stochastic Nonlinear Schr\"{o}dinger system]
{Global Solution and Blow-up of the Stochastic Nonlinear Schr\"{o}dinger system}

\author{Qi Zhang$^1$, Jinqiao Duan$^1$, Yong Chen$^2$}

\address{1. Department of Applied Mathematics, Illinois Institute of Technology, Chicago, IL 60616, USA}

\address{2. School of Science, Zhejiang Sci-Tech University, Hangzhou 310018, PR China}

\thanks{E-mail: Zhangqi1516@163.com, \ duan@iit.edu, \ youngchen329@126.com}

\maketitle

\begin{abstract}
We study the stochastic Nonlinear Schr\"{o}dinger system with multiplicative white noise in energy space $H^1$. Based on deterministic and stochastic Strichartz estimates, we prove the local well-posedness and uniqueness of mild solution. Then we prove the global well-posedness in the mass subcritical case and the defocusing case. For the mass subcritical case, we also investigate the global existence when the $L^2$ norm of initial value is small enough. In addition, we study the blow-up phenomenon and give a sharp criteria.
\end{abstract}

\vskip 0.2truein

\noindent {\it Keywords:}  Nonlinear Schr\"{o}dinger system; Well-posedness; Stochastic partial differential equations; Blow-up.

\vskip 0.2truein

\section{Introduction and Results}
We consider the well-posedness and the blow-up phenomenon of the following stochastic Nonlinear Schr\"{o}dinger system
\begin{equation}\label{CNLE}
  \left\{
   \begin{aligned}
   & idu + (\Delta u+ (\lambda_{11}|u|^{2\sigma}+\lambda_{12}|v|^{\sigma+1}|u|^{\sigma-1})u)dt=u\circ \phi_1 dW(t),\\
   & idv + (\Delta v+ (\lambda_{21}|v|^{\sigma -1}|u|^{\sigma+1}+\lambda_{22}|u|^{2\sigma})v)dt=v\circ \phi_2 dW(t),\\
   & u(0,x)=u_0(x), v(0,x)=v_0(x),
   \end{aligned}
   \right.
\end{equation}
where the coefficients $\lambda_{ij} \in \mathbb{R}$ for $i,j=1,2$, $W(t)$ is a cylindrical Wiener process in $L^2 \left( \mathbb{R}^N \right)$ with filtration $\left( \mathcal{F}_t \right)_{t\geqslant 0}$, the notation $\circ$ stands for Stratonovitch integral, and $\phi_1,\phi_2$ are Hilbert-Schmidt operators from $L^2(\mathbb{R}^N)$ into $H^1(\mathbb{R}^N)$.

When $\phi_1,\phi_2=0$, (\ref{CNLE}) reduces to the deterministic nonlinear Schr\"{o}dinger system
\begin{equation}\label{DCNLE}
  \left\{
   \begin{aligned}
   & i\partial_t u = - \Delta u -(\lambda_{11}|u|^{2\sigma}+\lambda_{12}|v|^{\sigma+1}|u|^{\sigma-1})u,\\
   & i\partial_t v = - \Delta v -(\lambda_{21}|v|^{\sigma -1}|u|^{\sigma+1}+\lambda_{22}|u|^{2\sigma})v,\\
   & u(0,x)=u_0(x), v(0,x)=v_0(x).
   \end{aligned}
   \right.
\end{equation}
In physics, the nonlinear Schr\"{o}dinger system (\ref{DCNLE}) is a important model and appears in many branches of physics, especially in Bose-Einstein condensation (BEC). In studying the BEC, the nonlinear Schr\"{o}dinger system can describe the propagation of wave function with interaction between two component, or the spin and motion of a particle (see \cite{BAOzw,MBBECPRL}). It is also an important model in nonlinear optics (see \cite{AA}). The coefficients $\lambda_{ij}\in \mathbb{R}$ are interaction constants between the $i-$th and $j-$th component in the system with $\lambda_{ij}=\lambda_{ji}$ ($\lambda_{ij}<0$ for defocusing case and $>0$ for focusing case). In physically, the blow-up phenomenon represents the wave function collapse, see \cite{SBlowup}. As a mathematical problem, The deterministic nonlinear Schr\"{o}dinger system (\ref{DCNLE}) has been extensively studied in the last decade. Many papers work on the deterministic nonlinear Schr\"{o}dinger system in different aspects. Results about the solitary wave solutions and ground state have been studied extensively, see \cite{GZZJDE, LWNcoupled}. We refer to Refs. \cite{CGCNLS, FMsharpGN, MZCNLS} for well-posedness and finite time blow-up for Cauchy problem.

In many circumstances, spatial and temporal fluctuations of the parameters of the medium have to be taken into account. It often occurs through a random potential, or describes the propagation of dispersive waves in nonhomogeneous or random media. In these cases, the stochastic Schr\"{o}dinger equations are introduced. For physical interpretations, we refer to \cite{BC,BG} and references therein.  From the mathematical point of view, the stochastic Schor\"{o}dinger equations is also an important problem. There are many well known results for the single stochastic Schor\"{o}dinger equation. See e.g. \cite{BDL2,BDH1} for the global well-posedness, \cite{BDB} for noise effects on blow-up, \cite{CHartree} for Hartree type nonlinear term. In \cite{BRZ1,BRZ2}, the authors consider the well-posedness problem of Stochastic nonlinear Schr\"{o}dinger equations with linear multiplicative
noise by rescaling approach. In \cite{BSSH}, Brze$\acute{z}$niak and Millet considered the  global existence and uniqueness of stochastic Schor\"{o}dinger equation on a compact Riemannian manifold based on a new Strichartz estimate for the stochastic convolution. The Schr\"{o}dinger equations stochastic with rougher noise already been considered by Oh, Pocovnicu, and Wang in \cite{FOW,OPWSNLS}.

Motivated by these problems, we are interested in the nonlinear Schr\"{o}dinger system with multiplicative cylindrical Wiener noise $W(t)$. The noise $W(t)$ acting as a random potential that are  dependent on $t$ and $x$. The aim of this paper is to state analogous results for the Cauchy problem of stochastic coupled system (\ref{CNLE}), including the local and global well-posedness and the blow-up phenomenon of $H^1$ solution.  As we know, it is the first attempt to study the stochastic coupled nonlinear Schr\"{o}dinger system. In our proof of the local well-posedness, we employ stochastic Strichartz estimates from \cite{BSSH,HSSE} to weaken the assumption of $\phi_1,\phi_2$ from \cite{BDL2,BDH1}. Moreover, thank to stochastic Strichartz estimates, we use the admissible pair $(r,2\sigma+2)$ directly, and adapt the deterministic fixed point argument in a ball of $C([0,T];H^1(\mathbb{R}^N))\cap L^r([0,T];W^{1,2\sigma+2}(\mathbb{R}^N))$. Then by the conversation of mass and the estimate of energy, we prove the global existence for the mass subcritical case and the defocusing case. For mass critical case, thanks to the sharp Gagliardo-Nirenberg inequality, we show the global existence when $\sqrt{\lambda_{11}} \|u_0\|^2_{L^2}+\sqrt{\lambda_{22}} \|v_0\|^2_{L^2}$ is sufficient small almost surely. The well-posedness results are similar with deterministic case \cite{CGCNLS, FMsharpGN, MZCNLS}. Moreover, the sharp criteria for blow-up is also discussed via generalizing the variance identity.

\subsection{Notations and Preliminaries}
Throughout the paper, the following notations and assumptions will be used. We assume that $0 \leqslant \sigma < \frac{2}{(N-2)^{+}}$ ($\frac{2}{(N-2)^{+}}=\infty$ when $N=1,2$, and $\frac{2}{(N-2)^{+}}=\frac{2}{N-2}$ when $N\geqslant3$). Let $\Lambda = \left(\begin{array}{l}{\lambda_{11} \  \lambda_{12}} \\ {\lambda_{21} \ \lambda_{22}}\end{array}\right)$ be the coefficient matrix in the following.

For $p \in \mathbb{N}^{\ast}$, $L^p$ is the Lebesgue space of complex valued functions. Moreover, $W^{1,p},H^1$ denotes the usual Sobolev space $W^{1,p}(\mathbb{R}^N)$ and $H^1\left(\mathbb{R}^N\right)$. The space $L^2$ is endowed with the inner product $(u, v)_{L^2}=\operatorname{Re} \int_{\mathbb{R}^{N}} u(x) \overline{v}(x) dx$. We also define the pseudo-conformal space $\Sigma=\left\{u \in H^{1}\left(\mathbb{R}^{N}\right) :|\cdot| u(\cdot) \in L^{2}\left(\mathbb{R}^{N}\right)\right\}
$ with norm $\|u\|_{\Sigma}^{2}=\|u\|_{H^{1}}^{2}+\|x u\|_{L^{2}}^{2}$.

Equation (\ref{CNLE}) has the following quantities: mass $M$, Hamiltonian $H$, variance $V$, and momentum $G$,
\begin{equation*}
M(u,v)= \int_{\mathbb{R}^N}|u|^2+|v|^2 dx,
\end{equation*}
\begin{equation*}
  \left.
   \begin{aligned}
   H(u,v) =&\frac{1}{2}\int_{\mathbb{R}^N} \left( |\nabla u|^2+|\nabla v|^2 \right) dx-\frac{1}{2+2\sigma} \int_{\mathbb{R}^N} \left( \lambda_{11} |u|^{2+2\sigma} \right. \\
       &\left. +\lambda_{22}|v|^{2+2\sigma}+2\lambda_{12}|v|^{\sigma+1}|u|^{\sigma+1} \right) dx
   \end{aligned}
   \right.
\end{equation*}
\begin{equation*}
V(u,v)=\int_{\mathbb{R}^{N}}|x|^{2} \left( |u(x)|^{2}+|v(x)|^{2} \right) d x, \quad u,v \in \Sigma .
\end{equation*}
\begin{equation*}
G(u,v)=\operatorname{Im} \int_{\mathbb{R}^{N}} u(x) x \cdot \nabla \overline{u}(x)+v(x) x \cdot \nabla \overline{v}(x) d x, \quad u, v \in \Sigma .
\end{equation*}

Given two separable Hilbert Spaces $H$ and $\tilde{H}$, and let $\phi: H \rightarrow \tilde{H}$ be a bounded linear operator. The operator $\phi$ is called Hilbert-Schmidt operator if there is an orthonormal basis $\left(e_{k}\right)_{k \in \mathbb{N}}$ in $H$ such that
\begin{equation*}
\operatorname{tr}\left(\phi^{*} \phi\right)=\sum_{k \in \mathbb{N}}\|\phi e_{k}\|_{\tilde{H}}^{2} < \infty.
\end{equation*}
If $\phi$ is a Hilbert-Schmidt operator, then
\begin{equation*}
\|\phi\|_{HS(H, \tilde{H})}=\left( \sum_{k \in \mathbb{N}}\|\phi e_{k}\|_{\tilde{H}}^{2} \right)^{\frac{1}{2}}
\end{equation*}
is called the Hilbert-Schmidt norm of $\phi$.

The cylindrical Wiener process on $L^2(\mathbb{R}^N)$ is defined by $W(t):=\sum_{k=0}^{\infty} e_k(x) B_k (t)$, here $(B_k)_{k \geqslant 0}$ is a sequence of independent standard Brownian motions, and $(e_k)_{k \geqslant 0}$ is an orthonormal basis of $L^2(\mathbb{R}^N)$. Assume that $F_{\phi_i}(x):=\sum_{k=0}^{\infty}(\phi e_k(x))^2 \in W^{1,\infty}(\mathbb{R}^N),i=1,2 $. Then the Stratonovitch integral can be written as following It\^o from:
\begin{align*}
\int_0^t u(s,x) \circ \phi_i dW(s,x)= & \int_0^t u(s,x)\phi_i dW(s,x)-\frac{1}{2}\int_0^t u(s,x) F_{\phi_i}(x)dt \\
                    =& \sum_{k=0}^{\infty}\int_0^t u(s,x) \phi_i e_k(x) dB_k (s)- \frac{1}{2}\int_0^t u(s,x) F_{\phi_i}(x)ds, \quad i=1,2.
\end{align*}

The unitary group $S(t):=e^{it\Delta}$ on $H^1(\mathbb{R}^N)$ is given by
\begin{equation*}
e^{it\Delta}f(x)=\frac{1}{(4\pi it)^{\frac{N}{2}}}\int_{\mathbb{R}^N} e^{\frac{i|x-y|^2}{4t}}f(y)dy.
\end{equation*}
By using the unitary group $S(t):=e^{it\Delta}$, the equivalent integral form of (\ref{CNLE}) is given by
\begin{equation}\label{CNLEI}
  \left\{
   \begin{aligned}
   u(t)=&S(t)u_0 +i \int_0^t S(t-s) (\lambda_{11}|u|^{2\sigma}+\lambda_{12}|v|^{\sigma+1}|u|^{\sigma-1})uds\\
        &+i\int_0^t S(t-s)u(s)\phi_1dW(s)-\frac{i}{2}\int^t_0S(t-s)u(s)F_{\phi_1}ds,\\
   v(t)=&S(t)v_0 +i \int_0^t S(t-s) (\lambda_{21}|v|^{\sigma-1}|u|^{\sigma+1}+\lambda_{22}|u|^{2\sigma})vds\\
        &+i\int_0^t S(t-s)v(s)\phi_2 dW(s)-\frac{i}{2}\int^t_0S(t-s)v(s)F_{\phi_2}ds.
   \end{aligned}
   \right.
\end{equation}

In the proof of local well-posedness, the following deterministic and stochastic Strichartz estimates is important in applying the contraction mapping argument to prove the local well-posedness of (\ref{CNLEI}). We say that a pair $(r, p)$ is admissible if $\frac{2}{r}+\frac{N}{p}=\frac{N}{2}$, and $2 \leq p \leq \frac{2 N}{(N-2)^{+}}$.
Now we state the well-known deterministic Strichartz estimates:
\begin{lemma}(Deterministic Strichartz estimates)
(i)Let $(r, p)$ be an admissible pair. For every $\varphi \in L^{2}\left(\mathbb{R}^{N}\right),$ the function $ S(t)\varphi \in L^{r}\left(\mathbb{R}, W^{1,p}\left(\mathbb{R}^{N}\right)\right) \cap C\left(\mathbb{R}, H^{1}\left(\mathbb{R}^{N}\right)\right)$. Furthermore, there exists a constant $C$ such
that
\begin{equation*}
\|S(\cdot) \varphi\|_{L^{r}\left(\mathbb{R}, W^{1,p}\left(\mathbb{R}^{N}\right)\right)} \leq C\|\varphi\|_{H^{1}\left(\mathbb{R}^{N}\right)}.
\end{equation*}
(ii) Let $I$ be an interval of $\mathbb{R}$ (bounded or not), $J=\overline{I},$ and $t_{0} \in J .$ If $(\gamma, \beta)$ is an admissible pair, $\frac{1}{\gamma}+\frac{1}{\gamma^{\prime}}=\frac{1}{\rho}+\frac{1}{\beta^{\prime}}=1$ and $f \in L^{\beta^{\prime}}\left(I, L^{\beta^{\prime}}\left(\mathbb{R}^{N}\right)\right),$ then for every admissible pair $(r, p)$, the function $G_{f}(t)=\int_{t_{0}}^{t} S(t-s) f(s) ds \in L^{r}\left(I, W^{1,p}\left(\mathbb{R}^{N}\right)\right) \cap C\left(J, H^{1}\left(\mathbb{R}^{N}\right)\right)$. Furthermore, there exists a constant $C$ independent of $I$ such that
\begin{equation*}
\left\|G_{f}\right\|_{L^{r}\left(I, W^{1,p}\left(\mathbb{R}^{N}\right)\right)} \leq C\|f\|_{L^{\gamma^{\prime}} \left(I, W^{1,\beta^{\prime}}\left(\mathbb{R}^{N}\right)\right)}.
\end{equation*}
\end{lemma}
Moreover, in order to deal the stochastic term in equation (\ref{CNLEI}), we need the following stochastic Strichartic estimates, see Brze$\acute{z}$niak and Millet \cite{BSSH}, or Hornung \cite{HSSE,HSSE1}.

\begin{lemma}\label{SSE}
(Stochastic Strichartz estimates) Let $T_1>0$, $\rho >1$.
For every predictable processes $\Phi \in L^{\rho}(\Omega,L^2([0,T_1]; HS(L^2(\mathbb{R}^N));H^1(\mathbb{R}^N))))$, the process
\begin{equation*}
J_{\left[0, T_{1}\right]} \Phi(t):=\int_{0}^{t} S(t-s) \Phi(s) \mathrm{d} W(s), \quad t \in\left[0, T_{1}\right]
\end{equation*}
is continuous and $\left( \mathcal{F}_t \right)_{t\geqslant 0}$-adapted in $H^1(\mathbb{R}^N)$ and $\left( \mathcal{F}_t \right)_{t\geqslant 0}$-predictable in $W^{1,p}(\mathbb{R}^N)$. Moreover, for every admissible pair $(r, p)$, there exists a constant $C$ independent of $T_1$ such that
\begin{align*}
& \left\|J_{\left[0, T_{1}\right]} \Phi\right\|_{L^{\rho}\left(\Omega, L^{r}\left([0, T_{1}]; W^{1,p}(\mathbb{R}^N)\right)\right.}+\left\|J_{\left[0, T_{1}\right]}^{T_{0}} \Phi\right\|_{L^{\rho}\left(\Omega, C\left(\left[0, T_{1}\right]; H^1(\mathbb{R}^N)\right)\right.} \\
\leqslant & C \|\Phi\|_{L^{\rho}\left(\Omega, L^{2}\left([0, T_{1}], HS(L^2(\mathbb{R}^N));H^1(\mathbb{R}^N)\right)\right)}.
\end{align*}
\end{lemma}

In order to use the Strichartz estimates, we choose the admissible pair $(r,2\sigma +2)$ such that $\frac{2}{r}+\frac{N}{2\sigma +2}=\frac{N}{2}$, and set the following Banach spaces:
\begin{equation*}
\mathcal{X}_t:=C([0,t];H^1(\mathbb{R}^N))\cap L^{r}([0,t];W^{1,2\sigma +2}(\mathbb{R}^N)),
\end{equation*}
\begin{equation*}
\mathcal{Y}_t:=C([0,t];L^2(\mathbb{R}^N))\cap L^{r}([0,t];L^{2\sigma +2}(\mathbb{R}^N)),
\end{equation*}
We also need the following sharp Gagliardo-Nirenberg inequality with best constant(see \cite{FMsharpGN}). From this inequality, we show that in mass critical cases, the sharp threshold for the global existence now depends on the ground state solutions of the elliptic system associated with (\ref{CNLE}).
\begin{lemma}(Sharp Gagliardo-Nirenberg inequality).
Let $\sigma < \frac{4}{(N-2)^{+}}$, $ \beta \geqslant 0$. Then for every $ u,v \in H^{1}\left(\mathbb{R}^{N}\right)$, we have
\begin{align*}
   & \|u\|^{2 \sigma+2}_{L^{2 \sigma+2}}+2 \beta\|u v\|^{\sigma+1}_{L^{\sigma+1}}+\|v\|^{2 \sigma+2}_{L^{2 \sigma+2}} \\
\leqslant & K_{o p t}\left(\|u\|^{2}_{L^2}+\|v\|_{L^2}^{2}\right)^{\sigma+1-\frac{\sigma N}{2}}\left(\|\nabla u\|_{L^2}^{2}+\|\nabla v\|_{L^2}^{2}\right)^{\frac{\sigma N}{2}},
\end{align*}
here the best constant
\begin{equation*}
K_{o p t}=\frac{2(\sigma +1)}{(N \sigma)^{N \sigma / 2}(2 \sigma+2-N \sigma)^{1-N \sigma / 2}} \frac{1}{\left(\|P\|_{L^2}^{2}+\|Q\|_{L^2}^{2}\right)^{\sigma}},
\end{equation*}
and $(P,Q)$ is the unique ground state solution of the nonlinear elliptic system
\begin{equation*}
\left\{\begin{array}{l}{-\Delta P+P-\left(|P|^{2 \sigma}+\beta|P|^{\sigma-1}|Q|^{\sigma +1}\right) P=0} \\ {-\Delta Q+Q-\left(|Q|^{2 \sigma}+\beta|Q|^{\sigma-1}|P|^{\sigma+1}\right) Q=0.}\end{array}\right.
\end{equation*}
\end{lemma}

\subsection{Main Results}

Applying the contraction mapping theorem in a suitable space and the Picard iteration method, we have the following local well-possdness and uniqueness result.

\begin{theorem}\label{LocalE}
Let $\sigma \in [0,\frac{2}{N}) \cup (\frac{1}{2},\frac{2}{(N-2)^{+}}) $, $\rho > 1$. Then given any $(u_0,v_0) \in L^{\rho}(\Omega, H^1(\mathbb{R}^N))$, there exists a stopping time $\tau^{\ast}(u_0,v_0)$ such that the equation (\ref{CNLE}) has a unique solution $(u,v) $ of equation (\ref{CNLE}) such that $ (u,z) \in L^{\rho}(\Omega ;C([0,\tau];(H^1(\mathbb{R}^N))^2)) $ for any $\tau < \tau^{\ast}(u_0,v_0)$. Moreover, we have almost surely
\begin{equation*}
\tau^{\ast}(u_0,v_0)=+\infty \quad \text{or} \quad  \lim_{t\nearrow \tau^{\ast}(u_0,v_0)}\sup_{s \leqslant t}||(u(s),v(s))||_{(H^1)^2}=+\infty .
\end{equation*}
\end{theorem}

\begin{remark}
If $N\leqslant 3$, than $\sigma \in [0,\frac{2}{(N-2)^{+}})$. Thus the restriction $\sigma \in [0,\frac{2}{N}) \cup (\frac{1}{2},\frac{2}{(N-2)^{+}})$ is not too strict in physics. 
\end{remark}

By the conservation of mass and the estimate of energy, following global well-posedness is established in this paper.

\begin{theorem}\label{GlobalE}
Under the same assumptions as Theorem \ref{LocalE}, suppose that one of following conditions holds\\
(i)   $\sigma < \frac{2}{N}$(mass subcritical),\\
(ii)  $\Lambda$ is nonpositive(defocusing),\\
(iii) $\sigma = \frac{2}{N}$(mass critical), $\Lambda$ is positive, and
\begin{equation*}
\left( \sqrt{\lambda_{11}} \|u_0\|^2_{L^2}+\sqrt{\lambda_{22}} \|v_0\|^2_{L^2} \right)<\frac{2}{\max \lbrace \lambda_{11},\lambda_{22} \rbrace K_{opt}} \quad a.s..
\end{equation*}
Then the solution of $(\ref{CNLE})$ given by Theorem \ref{LocalE} is global, that is $\tau^{*}\left(u_{0}\right)=+\infty,$ a.s..
\end{theorem}

\begin{remark}
In fact, we can obtain the global well-posedness and uniqueness of $L^{2}$-solution of (\ref{CNLE}) directly for each $\sigma $ from the the local well-posedness of $H^{1}$-solution and the conservation of mass.
\end{remark}

The following result gives a blow-up criteria for critical or supercritical case, i.e. $\frac{2}{N} \leqslant \sigma < \frac{2}{(N-2)^{+}}$.

\begin{theorem}\label{BLOWUP}
Assume that $\Lambda$ is negative, $\frac{2}{N} \leqslant \sigma < \frac{2}{(N-2)^{+}}$. Let all the assumptions of Theorem \ref{LocalE} hold. Let $(u,v)$ be a solution of equation (\ref{CNLE}) given by Theorem \ref{LocalE} with $(u_0,v_0) \in L^{\rho}(\Omega, \Sigma^2) $. Then for any stopping time $ \tau < \tau^{\ast}(u_0,v_0) $, the solution $(u,v) \in C \left( [0,\tau];\Sigma^2 \right)$ a.s.   Moreover, if
\begin{align}\label{Negative}
\mathbb{E}\left(V\left(u_0,v_0\right)\right)+ & 4 \mathbb{E}\left(G\left(u_0,v_0\right)\right) \overline{t}+8 \mathbb{E}\left(H\left(u_0,v_0\right)\right) \overline{t}^{2}\nonumber \\
  & \quad \quad + \frac{4}{3} \overline{t}^{3} \min_{i=1,2}\|F_{\phi_i}\|_{L^{\infty}} \mathbb{E}\left(M\left(u_0,v_0\right)\right)<0
\end{align}
for some $\overline{t}>0$, then $\mathbb{P}\left(\tau^{*}\left(u_{0},v_0\right) \leqslant \overline{t}\right)>0$.
\end{theorem}

\begin{remark}
Clearly, when $\|F_{\phi_1}\|_{L^{\infty}}$ or $\|F_{\phi_2}\|_{L^{\infty}}$ is small enough, the solution will blow up in finite time for initial data with negative energy almost surely. It is very similarly with deterministic case.
\end{remark}

This paper is organized as follows. In Section 2, we first address the proof of local well-posedness in subsection 2.1, then we prove the global well-posedness in subsection 2.2. In section 3, we shall prove the blow-up result.

\section{Well-posedness}
\subsection{Local Existence}
This subsection is mainly devoted to the proof of local existence $H^1$. We will define the solutions of (\ref{CNLE}) in Banach space $ L^{\rho}(\Omega; (\mathcal{X}_{T_0})^2)$ with $\rho \geqslant 1$. Let $E_T:=\lbrace (u,v)\in L^{\rho}(\Omega; (\mathcal{X}_{T})^2):(u(0),v(0))=(u_0,v_0)), \|u\|_{L^{\rho}(\Omega;\mathcal{X}_T)}+\|v\|_{L^{\rho}(\Omega;\mathcal{X}_T)} < M \rbrace $ be a subset of $ L^{\rho}(\Omega; (\mathcal{X}_{T})^2)$ with $M=2C \mathbb{E}(\|u_0\|^{\rho}_{H^1}+\|u_0\|^{\rho}_{H^1})$. For the nonlinear index $\sigma$, we can divide the situation into two cases, case $\sigma \in (\frac{1}{2},\frac{2}{(N-2)^{+}})$, and case $\sigma \in [0,\frac{2}{N}) \setminus (\frac{1}{2}, \frac{2}{(N-2)^{+}})$. For case $\sigma \in (\frac{1}{2},\frac{2}{(N-2)^{+}})$, $E_T$ equipped with the distance
\begin{equation*}
d\left( (u_1,v_1), (u_2,v_2) \right)=\|u\|_{L^{\rho}(\Omega;\mathcal{X}_T)}+\|v\|_{L^{\rho}(\Omega;\mathcal{X}_T)}.
\end{equation*}
here $(E_T,d)$ is a a complete metric space.
For case $\sigma \in [0,\frac{2}{N}) \setminus (\frac{1}{2}, \frac{2}{(N-2)^{+}})$, $E_T$ equipped with the distance
\begin{equation*}
d\left( (u_1,v_1), (u_2,v_2) \right)=\|u\|_{L^{\rho}(\Omega;\mathcal{Y}_T)}+\|v\|_{L^{\rho}(\Omega;\mathcal{Y}_T)},
\end{equation*}
Then from Lemma 3.19 in \cite{HSSE1}, $(E_T,d)$ is a complete metric space.
We then prove the local existence by using the contraction mapping argument pathwisely in metric space $(E_T,d)$.

In order to handle the nonlinear term in the concentration mapping argument, we need to truncate the nonlinear term in (\ref{CNLEI}), and consider the mild solution of following truncated equation:
\begin{equation}\label{CNLETI}
  \left\{
   \begin{aligned}
   u(t)=&S(t)u_0 +i \int_0^t S(t-s)\theta_R(u,v)(\lambda_{11}|u|^{2\sigma}+\lambda_{12}|v|^{\sigma+1}|u|^{\sigma-1})uds\\
        &+i\int_0^t S(t-s)u(s)\phi_{1}dW(s)-\frac{i}{2}\int^t_0S(t-s)u(s)F_{\phi_{1}}ds,\\
   v(t)=&S(t)v_0 +i \int_0^t S(t-s)\theta_R(u,v)(\lambda_{21}|v|^{\sigma-1}|u|^{\sigma+1}+\lambda_{22}|v|^{2\sigma})vds\\
        &+i\int_0^t S(t-s)v(s)\phi_{2}dW(s)-\frac{i}{2}\int^t_0S(t-s)v(s)F_{\phi_{2}}ds.
   \end{aligned}
   \right.
\end{equation}
The truncated $\theta_R$ is defined by
\begin{equation*}
\theta_R (u,v)=
\left\{
\begin{array}{rcl}
\Phi_R(\|(u,v)\|_{\mathcal{X}_{T}}),      &   \text{if}   & {\sigma \in (\frac{1}{2}, \frac{2}{(N-2)^{+}}) }\\
\Phi_R(\|(u,v)\|_{\mathcal{Y}_{T}}),       &  \text{if}    & {\sigma \in [0, \frac{2}{N}) \setminus (\frac{1}{2}, \frac{2}{(N-2)^{+}}). }
\end{array} \right.
\end{equation*}
where $\Phi_R \in C^{\infty}_c(\mathbb{R})$ with $\operatorname{supp}  \Phi_R \subset (-2R,2R)$, $0\leqslant \Phi\leqslant1$ ,and
\begin{equation*}
\Phi_R (x)=
\left\{
\begin{array}{rcl}
1,      &   \text{if}   & {|x| \leqslant R}\\
0,       &  \text{if}    & {|x| \geqslant 2R.}
\end{array} \right.
\end{equation*}
For $(u,v)\in L^{\rho}(\Omega; (\mathcal{X}_{T_0})^2)$, we define the mapping $\mathcal{T}(u,v)$ as the right hand side of (\ref{CNLETI}) on $(E,d)$.

\begin{lemma}
Let $\phi_1, \phi_{2} \in HS(L^2;H^1)$. Then for each each $(u,v)\in (E_T,d)$, $\mathcal{T}(u,v)\in (E_T,d)$ provided that $T>0$ is chosen small enough.
\end{lemma}
\noindent{\bf Proof.}
for each $(u,v)\in (E_T,d)$, $\mathcal{T}(u,v)\in (E_T,d)$. By deterministic Strichartz estimates, we have for almost surely $\omega$,
\begin{align}
\|\mathcal{T}(u,v)\|_{\left( \mathcal{X}_{T}\right)^2} \lesssim & \|u_0\|_{H^1}+ \|v_0\|_{H^1} \\ \nonumber
 & + \| \theta_R(u,v)(\lambda_{11}|u|^{2\sigma}+\lambda_{12}|v|^{\sigma+1}|u|^{\sigma-1})u\|_{L^{r^{\prime}}\left([0, T] ; W^{1,\frac{2\sigma +2}{2\sigma +1}}\right)}\\ \nonumber
  & + \| \theta_R(u,v)(\lambda_{21}|v|^{\sigma+1}|u|^{\sigma-1}+\lambda_{22}|v|^{2\sigma})v \|_{L^{r^{\prime}}\left([0, T] ; W^{1,\frac{2\sigma +2}{2\sigma +1}}\right)}\\ \nonumber
  & + \left\|\int_{0}^{t} S(t-s)v \phi_{2}d W(s)\right\|_{\mathcal{X}_{T}} + \left\|\int_{0}^{t} S(t-s)u \phi_{1}d W(s)\right\|_{\mathcal{X}_{T}} \\ \nonumber
  & + \|u F_{\phi_1}\|_{L^{r^{\prime}}\left([0, T] ; W^{1,\frac{2\sigma +2}{2\sigma +1}}\right)}+ \| vF_{\phi_2}\|_{L^{r^{\prime}}\left([0, T] ; W^{1,\frac{2\sigma +2}{2\sigma +1}} \right)} \\ \nonumber
:= & \|u_0\|_{H^1}+\|v_0\|_{H^1}+T_1+T_2+T_3+T_4.
\end{align}
We define a stopping time by
\begin{equation*}
t^R := \inf \lbrace 0< t \leqslant T: \| (u,v) \|_{(\mathcal{X}_t)^2} > 2R \rbrace.
\end{equation*}
Since $\|(u,v)\|_{(\mathcal{Y}_t)^2} \leqslant \|(u,v)\|_{(\mathcal{X}_t)^2}$, then for each $\sigma\in [0, \frac{2}{N}] \cup [\frac{1}{2}, \frac{2}{(N-2)^{+}})$, we have $\theta_R(u,v)=0$ when $t> t^R$. For $T_1$, from H\"{o}lder inequality, we get
\begin{align}\label{T1}
T_1 \leqslant & \| (\lambda_{11}|u|^{2\sigma}+\lambda_{12}|v|^{\sigma+1}|u|^{\sigma-1})u\|_{L^{r^{\prime}}\left([0, t^R] ; W^{1,\frac{2\sigma +2}{2\sigma +1}}\right)}\\  \nonumber
 \lesssim &   \| u \|^{2\sigma}_{L^{\eta}\left([0, t^R] ; W^{1,2\sigma +2}\right)}\| u \|_{L^{r}\left([0,t^R];W^{1,2\sigma +2}\right)}  \\ \nonumber
  & \quad \quad +\| v \|^{\sigma}_{L^{\eta}\left([0, t^R] ; W^{1,2\sigma +2}\right)}\| u \|^{\sigma}_{L^{\eta}\left([0, t^R] ; W^{1,2\sigma +2}\right)} \| v_1 \|_{L^{r}\left([0, t^R];W^{1,2\sigma +2}\right)} \\ \nonumber
 \lesssim & T^{\frac{2\sigma}{\eta}} \left(\| u \|^{2\sigma}_{L^{\infty}\left([0, t^R] ; W^{1,2\sigma +2}\right)}+ \| v \|^{2\sigma}_{L^{\infty}\left([0, t^R] ; W^{1,2\sigma +2}\right)} \right) \| u \|_{\mathcal{X}_{t^R}}\\ \nonumber
 & \quad \quad +T^{\frac{2\sigma}{\eta}}\| v \|^{\sigma}_{L^{\infty}\left([0, t^R] ; W^{1,2\sigma +2}\left(\mathbb{R}^N \right)\right)}\| u \|^{\sigma}_{L^{\infty}\left([0, t^R] ; W^{1,2\sigma +2}\right)} \| v \|_{\mathcal{X}_{t^R}} \\ \nonumber
   \lesssim & T^{\frac{2\sigma}{\eta}} \left( \| u \|^{2\sigma}_{\mathcal{X}_{t^R}}+ \| v \|^{2\sigma}_{\mathcal{X}_{t^R}} \right) \| u \|_{\mathcal{X}_{t^R}} \\ \nonumber
   \lesssim & T^{\frac{2\sigma}{\eta}} R^{2\sigma +1} .
\end{align}
Similarly, for $T_2$ we have
\begin{equation}\label{T2}
T_2 \lesssim T^{\frac{2\sigma}{\eta}} R^{2\sigma +1}.
\end{equation}
For stochastic term $T_3$, since $F_{\phi_1} \in W^{1,\infty}$, $\phi_1 \in HS(L^2,H^1)$, for each $g \in L^2$, $u(t)\phi_1 g \in H^1$. Moreover, for each orthonormal basis $ (e_k)_{k \geqslant 0} \subseteq L^2(\mathbb{R}^N)$, by H\"{o}lder inequality, for each $t \in [0,T]$,
\begin{equation*}
\|u(t)\phi_1\|^2_{HS}=\sum_{k \geqslant 0} \| u(t)\phi_1 e_k \|^2_{H^1} \leqslant  \|u(t)\|^2_{H^1} \| F_{\phi_1} \|_{W^{1,\infty}}<\infty.
\end{equation*}
Thus $S u(s)\phi_1 \in L^2([0,T];HS(L^2,H^1))$. Then stochastic Strichartiz estimates yields that
\begin{align*}
\left( \mathbb{E} \left\|\int_{0}^{t} S(t-s)u \phi_{1}d W(s)\right\|^{\rho}_{\mathcal{X}_{T_0}} \right)^{\frac{1}{\rho}} \lesssim &  \|S(t-s)u(s)\phi_1\|_{L^{\rho}\left(\Omega, L^{2}\left([0, T], HS(L^2;H^1 \right)\right)}\\ \nonumber
 \leqslant & T^{\frac{1}{2}}\|F_{\phi_1}\|^{\frac{1}{2}}_{W^{1,\infty}}\|u\|_{L^{\rho}\left(\Omega, L^{\infty}\left([0, T];H^1 \right)\right)} \\ \nonumber
  \leqslant & T^{\frac{1}{2}}\|u\|_{L^{\rho}\left(\Omega, \mathcal{X}_T\right)}
\end{align*}
Similarly,
\begin{equation*}
\left( \mathbb{E}\left\|\int_{0}^{t} S(t-s)v \phi_{2}d W(s)\right\|^{\rho}_{\mathcal{X}_{T_0}} \right)^{\frac{1}{\rho}} \lesssim T^{\frac{1}{2}}\|v\|_{L^{\rho}\left(\Omega, \mathcal{X}_T\right)},
\end{equation*}
So we conclude that
\begin{equation}\label{T3}
\left( \mathbb{E}T_3^{\rho} \right)^{\frac{1}{\rho}} \lesssim T^{\frac{1}{2}} \left( \|u\|_{L^{\rho}(\Omega;\mathcal{X}_T)}+\| v \|_{L^{\rho}(\Omega;\mathcal{X}_T)} \right) .
\end{equation}
Since $F_{\phi_1}, F_{\phi_2}\in W^{1,\infty}$, by H\"{o}lder inequality, we have
\begin{align}\label{T4}
T_4 \lesssim & T^{1-\frac{2}{r}} \left( \|u\|_{L^r([0,T];W^{1,2+2\sigma})} \|F_{\phi_1}\|_{W^{1,\infty}} + \|v\|_{L^r([0,T];W^{1,2+2\sigma})} \|F_{\phi_2}\|_{W^{1,\infty}} \right) \\ \nonumber
    \lesssim & T^{1-\frac{2}{r}} \left( \|F_{\phi_1}\|_{W^{1,\infty}} + \|F_{\phi_2}\|_{W^{1,\infty}} \right) \left( \| u \|_{\mathcal{X}_{T}}+\| v\|_{\mathcal{X}_{T}} \right).
\end{align}
From (\ref{T1})-(\ref{T4}), we conclude that
\begin{align}\label{ET}
& \|\mathcal{T}(u,v)\|_{L^{\rho}(\left(\Omega, \mathcal{X}_{T}\right)^2)}  \\ \nonumber
\lesssim & \|u_0\|_{H^1}+\|v_0\|_{H^1}+ T^{\frac{2\sigma}{\eta}}R^{2\sigma+1} + \left( T^{1-\frac{2}{r}} + T^{\frac{1}{2}} \right) \left( \|u\|_{L^{\rho}(\Omega;\mathcal{X}_T)}+\| v \|_{L^{\rho}(\Omega;\mathcal{X}_T)} \right).
\end{align}
Therefore when $T>0$ is chosen small enough, $\mathcal{T}(u,v) \in (E_T, d)$ for each $(u,v)\in (E_T, d)$. $\square$

\begin{lemma}\label{PLocalTE}
Let $\phi_1, \phi_{2} \in HS(L^2(\mathbb{R}^N;\mathbb{R});H^1(\mathbb{R}^N))$. Then for each $\mathcal{F}_0$-measurable $(u_0,v_0)$ taking value in $(H^1(\mathbb{R}^N))^2$, and for any given $T_0 >0$, the equation (\ref{CNLETI}) has a unique $H^1$-solution $(u,v)$ staring from $(u_0,v_0)$ which is almost surely in $ (\mathcal{X}_{T_0})^2 $.
\end{lemma}

\noindent{\bf Proof.} We prove the local existence by showing that $\mathcal{T}$ is a contraction mapping on $(E_T,d)$, where $T>0$ will be chosen small enough.

\noindent{\bf Case 1: $  \sigma \in (\frac{1}{2}, \frac{2}{(N-2)^{+}})$}\\
Let $(u_1,v_1),(u_2,v_2) \in (E_T,d) $, then by deterministic Strichartz estimates, we get for any $0<T<T_0$ and almost surely $\omega$,
\begin{align}\label{ECNLTI}
   &\| \mathcal{T}(u_1,v_1)-\mathcal{T}(u_2,v_2)\|_{(\mathcal{X}_T)^2}\\ \nonumber
\lesssim & \| \theta_R(u_1,v_1)(\lambda_{11}|u_1|^{2\sigma}+\lambda_{12}|v_1|^{\sigma+1}|u_1|^{\sigma-1})u_1\\ \nonumber
  & \quad \quad -\theta_R(u_2,v_2)(\lambda_{11}|u_2|^{2\sigma}+\lambda_{12}|v_2|^{\sigma+1}|u_2|^{\sigma-1})u_2\|_{L^{r^{\prime}}\left([0, T] ; W^{1,\frac{2\sigma +2}{2\sigma +1}}\right)}\\ \nonumber
  &+ \|\theta_R(u_1,v_1)(\lambda_{21}|v_1|^{\sigma+1}|u_1|^{\sigma-1}+\lambda_{22}|v_1|^{2\sigma})v_1 \\ \nonumber
  & \quad \quad -\theta_R(u_2,v_2)(\lambda_{21}|v_2|^{\sigma+1}|u_2|^{\sigma-1}+\lambda_{22}|v_2|^{2\sigma})v_2 \|_{L^{r^{\prime}}\left([0, T] ; W^{1, \frac{2\sigma +2}{2\sigma +1}}\right)}\\ \nonumber
  &+ \left\|\int_{0}^{t} S(t-s)\left(u_{1}(s)-u_{2}(s)\right) \phi_{1}d W(s)\right\|_{L^r([0,T];W^{1,2\sigma +2}(\mathbb{R}^N))}\\ \nonumber
  &+ \left\|\int_{0}^{t} S(t-s)\left(v_{1}(s)-v_{2}(s)\right) \phi_{2}d W(s)\right\|_{L^r([0,T];W^{1,2\sigma +2}(\mathbb{R}^N))}\\ \nonumber
  &+ \|\left(u_{1}-u_{2}\right)F_{\phi_1}\|_{L^{r^{\prime}}\left([0, T] ; W^{1, \frac{2\sigma +2}{2\sigma +1}}\right)}+ \| \left(v_{1}-v_{2}\right) F_{\phi_2}\|_{L^{r^{\prime}}\left([0, T] ; W^{1, \frac{2\sigma +2}{2\sigma +1}}\right)}\\ \nonumber
  &:= A+B+C+D+E.
\end{align}
In order to estimate $A$, for $R>0$, we set
\begin{equation*}
t^R_i:=\inf \lbrace\ 0< t \leqslant T, \| (u_i,v_i) \|_{(\mathcal{X}_t)^2} > 2R \rbrace , \quad i=1,2.
\end{equation*}
Without loss of generality, we assume that $t^R_1 \leqslant t^R_2$. Then $[0,T]=[0,t^R_1]\cup[t^R_1,t^R_2]\cup[t^R_2,T]$, and
\begin{align*}
A \lesssim & \| (\theta_R(u_1,v_1)-\theta_R(u_2,v_2))\\
  &\quad \quad (\lambda_{11}|u_1|^{2\sigma}+\lambda_{12}|v_1|^{\sigma+1}|u_1|^{\sigma-1})u_1  \|_{L^{r^{\prime}}\left([0, T_1^R] ; W^{1, \frac{2\sigma +2}{2\sigma +1}}\right)} \\
  +& \| \theta_R(u_2,v_2)\left( (\lambda_{11}|u_1|^{2\sigma}+\lambda_{12}|v_1|^{\sigma+1}|u_1|^{\sigma-1})u_1 \right.\\
  & \quad \quad \left. -(\lambda_{11}|u_2|^{2\sigma}+\lambda_{12}|v_2|^{\sigma+1}|u_2|^{\sigma-1})\right) u_2 \|_{L^{r^{\prime}}\left([0, T_2^R] ; W^{1,\frac{2\sigma +2}{2\sigma +1}}\right)}\\
  +& \| \theta_R(u_2,v_2)^2)(\lambda_{11}|u_2|^{2\sigma}+\lambda_{12}|v_2|^{\sigma+1}|u_2|^{\sigma-1})u_2\|_{L^{r^{\prime}}\left([T_1^R, T_2^R] ; W^{1,\frac{2\sigma +2}{2\sigma +1}}\right)}\\
  :=& A_1+A_2+A_3.
\end{align*}
In order to estimate $A_1$, by Lemma 3.3 in \cite{BDL2}, we have
\begin{equation*}
A_1 \lesssim \| (u_1,v_1)-(u_2,v_2) \|_{(\mathcal{X}_{T})^2} \| (\lambda_{11}|u_1|^{2\sigma}+\lambda_{12}|v_1|^{\sigma+1}|u_1|^{\sigma-1})u_1  \|_{L^{r^{\prime}}\left([0, T_1^R] ; W^{1, \frac{2\sigma +2}{2\sigma +1}}\right)}.
\end{equation*}
For the second factor, by H\"{o}lder's inequality with $\frac{2\sigma+1}{2\sigma+2}=\frac{2 \sigma}{2\sigma+2}+\frac{1}{2\sigma+2}$, $\frac{1}{\gamma^{\prime}}=\frac{2 \sigma}{\eta}+\frac{1}{r}$ , we have
\begin{align*}
 & \| (\lambda_{11}|u_1|^{2\sigma}+\lambda_{12}|v_1|^{\sigma+1}|u_1|^{\sigma-1})u_1  \|_{L^{r^{\prime}}\left([0, T_1^R] ; W^{1,\frac{2\sigma +2}{2\sigma +1}}\right)}\\
 \lesssim & \left\|\left(\left|u_{1}\right|^{2 \sigma}+\left|v_{1}\right|^{\sigma+1}\left|u_{1}\right|^{\sigma-1}\right) u_{1}\right\|_{L^{\gamma^{\prime}}\left(\left[0, t_{1}^{R}\right] ; L^{\frac{2\sigma+2}{2\sigma+1}}\right)} \\
 & \quad \quad + \left\|\left(\left|u_{1}\right|^{2 \sigma}+\left|v_{1}\right|^{\sigma+1}\left|u_{1}\right|^{\sigma-1}\right) \nabla u_{1}\|_{L^{\gamma}\left(\left[0, t_{1}^{R}\right] ; L^{s^{\prime}}\right)}  + \|\left|v_{1}\right|^{\sigma}\left|u_{1}\right|^{\sigma} \nabla v_{1}\right\|_{L^{\gamma}\left(\left[0, t_{1}^{R}\right] ; L^{\frac{2\sigma+2}{2\sigma+1}}\right)} \\
 \lesssim & \left(\left\|u_{1}\right\|_{L^{\eta}\left(\left[0, t_{1}^{R}\right] ; L^{2\sigma+2}\right)}^{2 \sigma}+\left\|v_{1}\right\|_{L^{\eta}\left(\left[0, t_{1}^{R}\right] ; L^{2\sigma+2}\right)}^{\sigma+1}\left\|u_{1}\right\|_{L^{\eta}\left(\left[0, t_{1}^{R}\right];L^{2\sigma+2}\right)}^{\sigma-1}\right)\left\|u_{1}\right\|_{L^{r}\left(\left[0, t_{1}^{R}\right] ; W^{1, 2\sigma +2}\right)} \\
 & \quad \quad + \left\|v_{1}\right\|_{L^{\eta}\left(\left[0, t_{1}^{R}\right] ; L^{2\sigma+2}\right)}^{\sigma}\left\|u_{1}\right\|_{L^{\eta}\left(\left[0, t_{1}^{R}\right] ; L^{2\sigma+2}\right)}^{\sigma}\left\|v_{1}\right\|_{L^{r}\left(\left[0, t_{1}^{R}\right] ; W^{1, 2\sigma+2}\right)} \\
 \lesssim & T^{\frac{2\sigma}{\eta}} \left(\left\|u_{1}\right\|^{2\sigma}_{L^{\infty}\left(\left[0, t_{1}^{R}\right] ; L^{2\sigma+2}\left(\mathbb{R}^{N}\right)\right)}+\left\|v_{1}\right\|^{2\sigma}_{L^{\infty}\left(\left[0, t_{1}^{R}\right] ; L^{2\sigma+2}\right)}\right)\left\|u_{1}\right\| \chi_{t_{1}^{R}} \\
  & \quad \quad +T^{\frac{2\sigma}{\eta}}\left\|v_{1}\right\|^{\sigma}_{L^{\infty}\left(\left[0, t_{1}^{R}\right] ; L^{2\sigma+2}\left(\mathbb{R}^{N}\right)\right)} \|u_{1}\|^{\sigma}_{L^{\infty}\left(\left[0, t_{1}^{R}\right] ; L^{2\sigma+2} \right) } \left\|v_{1}\right\|_{\mathcal{X}_{t_{1}^{R}}} \\
  \lesssim & T^{\frac{2\sigma}{\eta}} \left( \| u_1 \|^{2\sigma}_{{\mathcal{X}_{t^R_1}}}+ \| v_1 \|^{2\sigma}_{{\mathcal{X}_{t^R_1}}} \right) \left( \| u_1 \|_{\mathcal{X}_{t^R_1}}+\| v_1 \|_{\mathcal{X}_{t^R_1}} \right).
\end{align*}
where we have used the Sobolev embedding $H^1(\mathbb{R}^N)\hookrightarrow L^{2\sigma +2}(\mathbb{R}^N)$ since $2\sigma +2<\frac{2N}{N-2}$. Thus we deduce
\begin{equation}{\label{A1}}
  \left.
   \begin{aligned}
 A_1 \lesssim & \| (u_1,v_1)-(u_2,v_2) \|_{(\mathcal{X}_{T})^2} T^{\frac{2\sigma}{\eta}} \left( \| u_1 \|^{2\sigma}_{\mathcal{X}_{t^R_1}}+\| v_1 \|^{2\sigma}_{\mathcal{X}_{t^R_1}} \right) \left( \| u_1 \|_{\mathcal{X}_{t^R_1}}+\| v_1 \|_{\mathcal{X}_{t^R_1}} \right) \\
  \lesssim & T_0^{\frac{2\sigma}{\eta}}R^{2\sigma +1} \| (u_1,v_1)-(u_2,v_2) \|_{(\mathcal{X}_{T})^2}.
   \end{aligned}
   \right.
\end{equation}
For $A_3$, since $\theta_R(u_1,v_1)=0$ if $t \in \left(  t^R_1,t^R_2 \right) $, we use the same way for $A_1$ to estimate $A_3$ and get
\begin{equation}\label{A3}
 A_3 \lesssim T^{\frac{2\sigma}{\eta}}R^{2\sigma+1} \| (u_1,v_1)-(u_2,v_2) \|_{(\mathcal{X}_{T})^2}.
\end{equation}

For $A_2$, note that $\Phi_{R}$ is bounded, we have
\begin{align}\label{A2}
A_2 & \lesssim | \lambda_{11}| \| (|u_1|^{2\sigma} u_1 -|u_2|^{2\sigma} u_2)\|_{L^{r^{\prime}}\left([0, t^R_1] ; W^{1,\frac{2\sigma +2}{2\sigma +1}}\right)}\nonumber\\
  & \quad + |\lambda_{12}| \|(|v_1|^{\sigma+1}- |v_2|^{\sigma+1}) |u_1|^{\sigma-1}u_1\|_{L^{\gamma^{\prime}}\left([0, t^R_1] ; W^{1,\frac{2\sigma +2}{2\sigma +1}}\right)}\nonumber\\
  & \quad + | \lambda_{12}|\||v_2|^{\sigma+1}(|u_1|^{\sigma-1}u_1-|u_2|^{\sigma-1}u_2) \|_{L^{\gamma^{\prime}}\left([0, t^R_1] ; W^{1,\frac{2\sigma +2}{2\sigma +1}}\right)}\nonumber\\
  & \lesssim \| (|u_1|^{2\sigma}+|u_2|^{2\sigma}+|v_1|^{2\sigma}+|v_2|^{2\sigma}) \nonumber\\
  & \quad \quad  \times(|u_1-u_2|+|v_1-v_2|) \|_{L^{r^{\prime}}\left([0, t^R_1] ; W^{1, \frac{2\sigma +2}{2\sigma +1}}\right)} \nonumber\\
  & \lesssim \|(\left|u_{1}\right|^{2 \sigma}+\left|u_{2}\right|^{2 \sigma}+\left|v_{1}\right|^{2 \sigma}+\left|v_{2}\right|^{2 \sigma}) \nonumber \\
  & \quad \quad  \times (\nabla(|u_1-u_2|+|v_1-v_2|)+|u_1-u_2|+|v_1-v_2|) \|_{L^{r^{\prime}}\left([0, t^R_1] ; L^{\frac{2\sigma +2}{2\sigma +1}}\right)} \nonumber\\
  & + \| (|u_1|^{2\sigma-1}\nabla u_1 +|u_2|^{2\sigma-1}\nabla u_2 +|v_1|^{2\sigma-1}\nabla v_1 +|v_2|^{2\sigma-1}\nabla v_2) \nonumber\\
  & \quad \quad  \times(|u_1-u_2|+|v_1-v_2|) \|_{L^{r^{\prime}}\left([0, t^R_1] ; L^{\frac{2\sigma +2}{2\sigma +1}}\right)} \nonumber\\
  & \lesssim T^{\frac{2\sigma}{\eta}} \left(  \sum_{i=1,2}\|u_i\|^{2\sigma}_{L^{\infty}\left([0, t^R_1] ; L^{2\sigma +2}\right)}+\sum_{i=1,2}\|v_i\|^{2\sigma}_{L^{\infty}\left([0, t^R_1] ; L^{2\sigma +2}\right)} \right) \nonumber\\
  & \quad \quad  \times \left( \|u_1-u_2\|_{L^{r}([0,t^R_1];W^{1,2\sigma +2}} +\|v_1-v_2 \|_{L^{r}([0,t^R_1];W^{1,2\sigma +2}} \right) \nonumber\\
  & \quad \quad + T^{\frac{2\sigma}{\eta}} \left(  \sum_{i=1,2}\|u_i\|^{2\sigma-1}_{L^{\infty}\left([0, t^R_1] ; L^{2\sigma +2}\right)} \left\|\nabla u_{i}\right\|_{L^{r}\left(\left[0, t_{1}^{R}\right] ; L^{2 \sigma+2}\right) }  \right.  \nonumber\\
  & \quad \quad \left. +\sum_{i=1,2}\|v_i\|^{2\sigma-1}_{L^{\infty}\left([0, t^R_1] ; L^{2\sigma +2}\right)} \left\|\nabla v_{i}\right\|_{L^{r}\left(\left[0, t_{1}^{R}\right] ; L^{2\sigma+2}\right) } \right) \nonumber\\
  & \quad \quad \times\left(\left\|u_{1}-u_{2}\right\|_{L^{\infty}\left(\left[0, t_{1}^{R}\right] ; L^{2\sigma+2}\right)}+\left\|v_{1}-v_{2}\right\|_{L^{\infty}\left(\left[0, t_{1}^{R}\right] ; L^{2\sigma+2}\right)}\right) \nonumber \\
  &\lesssim T^{\frac{2\sigma}{\eta}}R^{2 \sigma +1}( \| u_1-u_2 \|_{\mathcal{X}_{T}}+\| v_1-v_2 \|_{\mathcal{X}_{T}} ).
   \end{align}
Collecting (\ref{A1})-(\ref{A2}) shows that
\begin{equation}\label{A}
 A \lesssim T^{\frac{2\sigma}{\eta}}R^{2 \sigma +1}( \| u_1-u_2 \|_{\mathcal{X}_{T}}+\| v_1-v_2 \|_{\mathcal{X}_{T}} ).
\end{equation}
Similarly, we have
\begin{equation}\label{B}
 B \lesssim T^{\frac{2\sigma}{\eta}}R^{2 \sigma +1}( \| u_1-u_2 \|_{\mathcal{X}_{T}}+\| v_1-v_2 \|_{\mathcal{X}_{T}} ).
\end{equation}
In stochastic term $C$ and $D$, $S(t-s)(u_1(s)-u_2(s))\phi_1, S(t-s)(v_1(s)-v_2(s))\phi_2 \in HS(L^2;H^1)$. Thus stochastic Strichartz estimates yields that
\begin{align}\label{CD}
& (\mathbb{E}C^{\rho})^{\frac{1}{\rho}}+(\mathbb{E}D^{\rho})^{\frac{1}{\rho}} \\ \nonumber
 = & \mathbb{E} \left( \left\|\int_{0}^{t} S(t-s)\left(u_{1}(s)-u_{2}(s)\right) \phi_{1}d W(s)\right\|^{\rho}_{L^r([0,T];W^{1,2\sigma +2})} \right)^{\frac{1}{\rho}} \\ \nonumber
   & + \mathbb{E}\left( \left\|\int_{0}^{t} S(t-s)\left(v_{1}(s)-v_{2}(s)\right) \phi_{2}d W(s)\right\|^{\rho}_{L^r([0,T];W^{1,2\sigma +2})}\right)^{\frac{1}{\rho}} \\ \nonumber
  \lesssim &  \|S(t-s)(u_1(s)-u_2(s))\phi_1\|_{L^{\rho}\left(\Omega, L^{2}\left([0, T], HS(L^2;H^1)\right)\right)}\\ \nonumber
   & +\|S(t-s)(v_1(s)-v_2(s))\phi_2\|_{L^{\rho}\left(\Omega, L^{2}\left([0, T], HS(L^2;H^1)\right)\right)}\\ \nonumber
  \lesssim &  T^{\frac{1}{2}} \left( \|F_{\phi_1}\|_{W^{1,\infty}}^{\frac{1}{2}}+\|F_{\phi_2}\|_{W^{1,\infty}}^{\frac{1}{2}} \right) \left( \|(u_1(s)-u_2(s))\|_{L^{\rho}\left(\Omega, L^{\infty}\left([0, T], H^1 \right)\right)} \right. \\ \nonumber
   & \left. +\|S(t-s)(v_1(s)-v_2(s))\phi_2\|_{L^{\rho}\left(\Omega, L^{\infty}\left([0, T], H^1 \right)\right)} \right) \\ \nonumber
  \lesssim & T^{\frac{1}{2}} \left( \|F_{\phi_1}\|_{W^{1,\infty}}^{\frac{1}{2}}+\|F_{\phi_2}\|_{W^{1,\infty}}^{\frac{1}{2}} \right) d((u_1,v_1),(u_2,v_2)).
\end{align}
Since $F_{\phi_1},F_{\phi_1} \in W^{1,\infty}(\mathbb{R}^N)$, from H\"{o}lder inequality, we have
\begin{equation}\label{E}
  \left.
   \begin{aligned}
E  \lesssim & T^{1-\frac{2}{r}} \left( \| u_1-u_2 \|_{L^{r}([0,t^R_1];L^{2 \sigma+2}(\mathbb{R}^N)}\|F_{\phi_1}\|_{W^{1,\infty}(\mathbb{R}^N)} \right. \\
  &  \left. +\| v_1-v_2 \|_{L^{r}([0,t^R_1];L^{2\sigma +2}(\mathbb{R}^N)} ) \|F_{\phi_2}\|_{W^{1, \infty}(\mathbb{R}^N)} \right),\\
  \lesssim & T^{1-\frac{2}{r}} \left( \|F_{\phi_1}\|_{W^{1,\infty}} + \|F_{\phi_2}\|_{W^{1,\infty}} \right) \left( \| u_1-u_2 \|_{\mathcal{X}_{T}}+\| v_1-v_2 \|_{\mathcal{X}_{T}} \right).
   \end{aligned}
   \right.
\end{equation}
Combining estimates (\ref{A})-(\ref{E}), we get that for $\nu=min(1-\frac{2}{r},\frac{2\sigma}{\eta},\frac{1}{2})$,
\begin{equation*}
d(\mathcal{T}(u_1,v_1)-\mathcal{T}(u_2,v_2))\leqslant C(R,\|\phi_1\|_{W^{1,\infty}},\|\phi_2\|_{W^{1,\infty}})T^{\nu}d(u_1-u_2, v_1-v_2 ).
\end{equation*}
Thus $\mathcal{T}$ is a contraction mapping in $(E,d)$ provided $T$ is chosen small enough. 

\noindent{\bf Case 2: $ \sigma \in [0, \frac{2}{N}) \setminus (\frac{1}{2}, \frac{2}{(N-2)^{+}})$}\\
 Let $(u_1,v_1),(u_2,v_2) \in (E_T,d) $, then by deterministic Strichartz estimates, we get for any $0<T<T_0$ and almost surely $\omega$,
\begin{align}\label{ECNLTIL2}
   &\| \mathcal{T}(u_1,v_1)-\mathcal{T}(u_2,v_2)\|_{(\mathcal{Y}_T)^2}\\ \nonumber
\lesssim & \| \Phi_R(\|(u_1,v_1)\|_{\mathcal{Y}_t})(\lambda_{11}|u_1|^{2\sigma}+\lambda_{12}|v_1|^{\sigma+1}|u_1|^{\sigma-1})u_1\\ \nonumber
  & \quad \quad -\Phi_R(\|(u_2,v_2)\|_{\mathcal{Y}_t})(\lambda_{11}|u_2|^{2\sigma}+\lambda_{12}|v_2|^{\sigma+1}|u_2|^{\sigma-1})u_2\|_{L^{r^{\prime}}\left([0, T] ; L^{\frac{2\sigma +2}{2\sigma +1}}\right)}\\ \nonumber
  &+ \|\Phi_R(\|(u_1,v_1)\|_{\mathcal{Y}_t})(\lambda_{21}|v_1|^{\sigma+1}|u_1|^{\sigma-1}+\lambda_{22}|v_1|^{2\sigma})v_1 \\ \nonumber
  & \quad \quad -\Phi_R(\|(u_2,v_2)\|_{\mathcal{Y}_t})(\lambda_{21}|v_2|^{\sigma+1}|u_2|^{\sigma-1}+\lambda_{22}|v_2|^{2\sigma})v_2 \|_{L^{r^{\prime}}\left([0, T] ; L^{ \frac{2\sigma +2}{2\sigma +1}}\right)}\\ \nonumber
  &+ \left\|\int_{0}^{t} S(t-s)\left(u_{1}(s)-u_{2}(s)\right) \phi_{1}d W(s)\right\|_{L^r([0,T];L^{2\sigma +2})}\\ \nonumber
  &+ \left\|\int_{0}^{t} S(t-s)\left(v_{1}(s)-v_{2}(s)\right) \phi_{2}d W(s)\right\|_{L^r([0,T];L^{2\sigma +2})}\\ \nonumber
  &+ \|\left(u_{1}-u_{2}\right)F_{\phi_1}\|_{L^{r^{\prime}}\left([0, T] ; L^{ \frac{2\sigma +2}{2\sigma +1}}\right)}+ \| \left(v_{1}-v_{2}\right) F_{\phi_2}\|_{L^{r^{\prime}}\left([0, T] ; L^{ \frac{2\sigma +2}{2\sigma +1}}\right)}\\ \nonumber
  &:= \overline{A}+\overline{B}+\overline{C}+\overline{D}+\overline{E}.
\end{align}
In order to estimate $\overline{A}$, for $R>0$, we set
\begin{equation*}
t^R_i:=\inf \lbrace\ 0< t \leqslant T, \| (u_i,v_i) \|_{(\mathcal{Y}_t)^2} > 2R \rbrace , \quad i=1,2.
\end{equation*}
Without loss of generality, we assume that $t^R_1 \leqslant t^R_2$. Then $[0,T]=[0,t^R_1]\cup[t^R_1,t^R_2]\cup[t^R_2,T]$, and
\begin{align*}
\overline{A} \lesssim & \| (\Phi_R(\|(u_1,v_1)\|_{\mathcal{Y}_{t}})-\Phi_R(\|(u_2,v_2) \|_{\mathcal{Y}_{t}}))\\
  &\quad \quad (\lambda_{11}|u_1|^{2\sigma}+\lambda_{12}|v_1|^{\sigma+1}|u_1|^{\sigma-1})u_1  \|_{L^{r^{\prime}}\left([0, T_1^R] ; L^{ \frac{2\sigma +2}{2\sigma +1}}\right)} \\
  +& \| \Phi_R(\|(u_2,v_2) \|_{(\mathcal{Y}_{t})^2})\left( (\lambda_{11}|u_1|^{2\sigma}+\lambda_{12}|v_1|^{\sigma+1}|u_1|^{\sigma-1})u_1 \right.\\
  & \quad \quad \left. -(\lambda_{11}|u_2|^{2\sigma}+\lambda_{12}|v_2|^{\sigma+1}|u_2|^{\sigma-1})\right) u_2 \|_{L^{r^{\prime}}\left([0, T_2^R] ; L^{\frac{2\sigma +2}{2\sigma +1}}\right)}\\
  +& \| \Phi_R(\|(u_2,v_2)\|_{(\mathcal{Y}_t})^2)(\lambda_{11}|u_2|^{2\sigma}+\lambda_{12}|v_2|^{\sigma+1}|u_2|^{\sigma-1})u_2\|_{L^{r^{\prime}}\left([T_1^R, T_2^R] ; L^{\frac{2\sigma +2}{2\sigma +1}}\right)}\\
  :=& \overline{A}_1+\overline{A}_2+\overline{A}_3.
\end{align*}
In order to estimate $\overline{A}_1$, by Lemma 3.3 in \cite{BDL2}, we have
\begin{equation*}
\overline{A}_1 \lesssim \| (u_1,v_1)-(u_2,v_2) \|_{(\mathcal{Y}_{T})^2} \| (\lambda_{11}|u_1|^{2\sigma}+\lambda_{12}|v_1|^{\sigma+1}|u_1|^{\sigma-1})u_1  \|_{L^{r^{\prime}}\left([0, T_1^R] ; L^{\frac{2\sigma +2}{2\sigma +1}}\right)}.
\end{equation*}
For the second factor, by H\"{o}lder's inequality with $\frac{1}{r^{\prime}}=\frac{2\sigma+1}{r}+(1-\frac{N\sigma}{2})$, we have
\begin{align*}
& \| (\lambda_{11}|u_1|^{2\sigma}+\lambda_{12}|v_1|^{\sigma+1}|u_1|^{\sigma-1})u_1  \|_{L^{r^{\prime}}\left([0, T_1^R] ; L^{\frac{2\sigma +2}{2\sigma +1}}\right)}\\
\lesssim & T^{1-\frac{N \sigma}{2}} \left( \| u_1 \|^{2\sigma+1}_{L^{r}\left([0, T_1^R] ; L^{2\sigma +2}\right)}+\| v_1 \|^{2\sigma}_{L^{r}\left([0, T_1^R] ; L^{2\sigma +2}\right)}\| u_1 \|_{L^{r}\left([0, T_1^R] ; L^{2\sigma +2}\right)} \right) \\
  \lesssim & T^{1-\frac{N \sigma}{2}} R^{2\sigma+1}.
\end{align*}
Thus we deduce
\begin{equation}{\label{A1L2}}
\overline{A}_1  \lesssim  T^{1-\frac{N \sigma}{2}}R^{2\sigma +1} \| (u_1,v_1)-(u_2,v_2) \|_{(\mathcal{Y}_{T})^2}.
\end{equation}
For $\overline{A}_3$, since $\Phi_R(\|(u_1,v_1)\|_{(\mathcal{X}_t)^2})=0$ if $t \in \left(  t^R_1,t^R_2 \right) $, we use the same way for $\overline{A}_1$ to estimate $\overline{A}_3$ and get
\begin{equation}\label{A3L2}
 \overline{A}_3 \lesssim T^{1-\frac{N \sigma}{2}}R^{2\sigma+1} \| (u_1,v_1)-(u_2,v_2) \|_{(\mathcal{Y}_{T})^2}.
\end{equation}

For $\overline{A}_2$, note that $\Phi_{R}$ is bounded, we have
\begin{align}\label{A2L2}
\overline{A}_2 & \lesssim | \lambda_{11}| \| (|u_1|^{2\sigma} u_1 -|u_2|^{2\sigma} u_2)\|_{L^{r^{\prime}}\left([0, t^R_1] ; L^{\frac{2\sigma +2}{2\sigma +1}}\right)}\nonumber\\
  & \quad + |\lambda_{12}| \|(|v_1|^{\sigma+1}- |v_2|^{\sigma+1}) |u_1|^{\sigma-1}u_1\|_{L^{r^{\prime}}\left([0, t^R_1] ; L^{\frac{2\sigma +2}{2\sigma +1}}\right)} \nonumber\\
  & \quad + | \lambda_{12}|||v_2|^{\sigma+1}(|u_1|^{\sigma-1}u_1-|u_2|^{\sigma-1}u_2)| \|_{L^{\gamma^{\prime}}\left([0, t^R_1] ; L^{\frac{2\sigma +2}{2\sigma +1}}\right)}\nonumber\\
  & \lesssim \| (|u_1|^{2\sigma}+|u_2|^{2\sigma}+|v_1|^{2\sigma}+|v_2|^{2\sigma}))\nonumber\\
  & \quad \quad  \times(|u_1-u_2|+|v_1-v_2|) \|_{L^{r^{\prime}}\left([0, t^R_1] ; L^{\frac{2\sigma +2}{2\sigma +1}}\right)} \nonumber\\
  & \lesssim T^{1-\frac{N \sigma}{2}} \left(  \sum_{i=1,2}\|u_i\|^{2\sigma}_{L^r\left([0, t^R_1] ; L^{2\sigma +2}\left(\mathbb{R}^N\right)\right)}+\sum_{i=1,2}\|v_i\|^{2\sigma}_{L^r\left([0, t^R_1] ; L^{2\sigma +2}\left(\mathbb{R}^N\right)\right)} \right) \nonumber\\
  & \quad \quad  \times \left( \|u_1-u_2\|_{L^{r}([0,t^R_1];L^{2\sigma +2}} +\|v_1-v_2 \|_{L^{r}([0,t^R_1];L^{2\sigma +2}} \right) \nonumber\\
  &\lesssim T^{1-\frac{N \sigma}{2}}R^{2 \sigma}( \| u_1-u_2 \|_{\mathcal{Y}_{T}}+\| v_1-v_2 \|_{\mathcal{Y}_{T}} ).
   \end{align}
Collecting (\ref{A1L2})-(\ref{A2L2}) shows that
\begin{equation}\label{AL2}
 \overline{A} \lesssim T^{1-\frac{N \sigma}{2}}R^{2 \sigma +1}( \| u_1-u_2 \|_{\mathcal{Y}_{T}}+\| v_1-v_2 \|_{\mathcal{Y}_{T}} ).
\end{equation}
Similarly, we have
\begin{equation}\label{BL2}
 \overline{B} \lesssim T^{1-\frac{N \sigma}{2}}R^{2 \sigma +1}( \| u_1-u_2 \|_{\mathcal{Y}_{T}}+\| v_1-v_2 \|_{\mathcal{Y}_{T}} ).
\end{equation}
In stochastic term $\overline{C}$ and $\overline{D}$, $S(t-s)(u_1(s)-u_2(s))\phi_1, S(t-s)(v_1(s)-v_2(s))\phi_2 \in HS(L^2;H^1)$. Thus similar with (\ref{CD}), we have
\begin{equation}\label{CDL2}
(\mathbb{E}\overline{C}^{\rho})^{\frac{1}{\rho}}+(\mathbb{E}\overline{D}^{\rho})^{\frac{1}{\rho}}  \lesssim  T^{\frac{1}{2}} \left( \|F_{\phi_1}\|_{L^{\infty}}^{\frac{1}{2}}+\|F_{\phi_2}\|_{L^{\infty}}^{\frac{1}{2}} \right) d((u_1,v_1),(u_2,v_2)).   
\end{equation}
Since $F_{\phi_1},F_{\phi_1} \in W^{1,\infty}$, from H\"{o}lder inequality, we have
\begin{equation}\label{EL2}
  \left.
   \begin{aligned}
\overline{E}  \lesssim & T^{1-\frac{2}{r}} \left( \| u_1-u_2 \|_{L^{r}([0,t^R_1];L^{2 \sigma+2}}\|F_{\phi_1}\|_{L^{\infty}} \right. \\
  &  \left. +\| v_1-v_2 \|_{L^{r}([0,t^R_1];L^{2\sigma +2}} ) \|F_{\phi_2}\|_{L^{\infty}} \right),\\
  \lesssim & T^{1-\frac{2}{r}} \left( \|F_{\phi_1}\|_{W^{1,\infty}} + \|F_{\phi_2}\|_{W^{1,\infty}} \right) \left( \| u_1-u_2 \|_{\mathcal{Y}_{T}}+\| v_1-v_2 \|_{\mathcal{Y}_{T}} \right).
   \end{aligned}
   \right.
\end{equation}
Combining estimates (\ref{AL2})-(\ref{EL2}), we get that for $\nu=min(1-\frac{2}{r},1-\frac{N\sigma}{2},\frac{1}{2})$,
\begin{equation*}
d(\mathcal{T}(u_1,v_1)-\mathcal{T}(u_2,v_2))\leqslant C(R,\|\phi_1\|_{W^{1,\infty}},\|\phi_2\|_{W^{1,\infty}})T^{\nu}d(u_1-u_2, v_1-v_2 ).
\end{equation*}
Thus $\mathcal{T}$ is a contraction mapping in $(E,d)$ provided $T$ is chosen small enough.

Hence the truncated equation (\ref{CNLETI}) has a unique solution $(u^R, v^R)$ in $L^{\rho}(\Omega, (\mathcal{X}_{T})^2)$. Moreover, the solution can easily be extended to the whole interval $[0, T_0]$ by continuing the solution to $[T,2T]$, $[2T,3T]$, $\ldots$, and so on. The proof of Lemma \ref{PLocalTE} is thus completed. $\square$

\noindent{\bf Proof of Theorem \ref{LocalE}.}Let
\begin{equation*}
\tau_R (\omega)=\left\{\begin{array}{ll}{\inf(t \in [0,T_0], \|(u,v)\|_{(\mathcal{X}_{T})^2} \leqslant R)} & {\text { if } \sigma \in (\frac{1}{2}, \frac{2}{(N-2)^{+}}) } \\ {\inf(t \in [0,T_0], \|(u,v)\|_{(\mathcal{Y}_{T})^2} \leqslant R)} & {\text { if } \sigma \in [0, \frac{2}{N}) \setminus (\frac{1}{2}, \frac{2}{(N-2)^{+}}) .}\end{array}\right.
\end{equation*}
We denote for $m \in \mathbb{N}$ by $(u^m, v^m)$ the unique global solution of (\ref{CNLETI}) with $R=m$. Now we show that $\tau_m$ is increasing with $m$ and the $H^1$-solution of equation (\ref{CNLE}) can be defined by $(u,v)=(u^m, v^m)$ on $[0,\tau_m]$.
We need following lemma.
\begin{lemma}\label{truncated}
$(u^m, v^m)=(u^{m+1}, v^{m+1})$ for each $t\in [0, min(\tau_m, \tau_{m+1})]$ for a.e. $\omega \in \Omega$.
\end{lemma}
\noindent{\bf Proof.}
Fix $m\in\mathbb{N}$ and $T$, and let $\tau=min(\tau_m,\tau_{m+1})$. If $\tau \in T_0$, we define $U_m$ as solution on $[\tau, T_0]$ of the equation
\begin{equation*}
idU_m+\Delta U_m dt=U_m \phi_1 dW-\frac{i}{2}U_m F_{\phi}dt
\end{equation*}
with $U_m(\tau)=u_m(\tau)$.
And we define $V_m$ as solution on $[\tau, T_0]$ of the equation
\begin{equation*}
idV_m+\Delta V_m dt=V_m \phi_2 dW-\frac{i}{2}V_m F_{\phi}dt
\end{equation*}
with $U_m(\tau)=u_m(\tau)$.

Now we set
\begin{equation}
( \tilde{u}_{m}(t),\tilde{v}_{m}(t) )=\left\{\begin{array}{ll}{(u_{m}(t),v_{m}(t))} & {\text { if } t \in[0, \tau]} \\ {(U_{m}(t),V_{m}(t))} & {\text { if } t \in\left[\tau, T_{0}\right]}\end{array}\right.
\end{equation}
For $m+1$, we can define $( \tilde{u}_{m+1}(t),\tilde{v}_{m+1}(t) )$by same way. Now we will show that when $T$ is small enough, $( \tilde{u}_{m}(t),\tilde{v}_{m}(t) )=( \tilde{u}_{m+1}(t),\tilde{v}_{m+1}(t) )$ on $[0,T]$, and the lemma will follow by a reiteration argument. For $t\in [0,T]$,
\begin{align*}
  & \tilde{u}_{m+1}(t)-\tilde{u}_{m}(t)\\
 =&-i \int_{0}^{t \wedge \tau} S(t-s)\left(\left(|\lambda_{11}| \left|\tilde{u}_{m+1}(s)\right|^{2\sigma}+|\lambda_{12}|\left|\tilde{v}_{m+1}(s)\right|^{\sigma+1}\left|\tilde{u}_{m-1}(s)\right|^{\sigma-1}\right) \tilde{u}_{m+1}(s) \right. \\
 & \left. -\left(\left| \lambda_{11}| \tilde{u}_{m}(s)\right|^{2\sigma}+|\lambda_{12}|\left|\tilde{v}_{m}(s)\right|^{\sigma+1}\left|\tilde{v}_{m}(s)\right|^{\sigma-1}\right) \tilde{u}_{m}(s)\right) ds  \\
 & -i \int_{0}^{t} S(t-s)\left(\left(\tilde{u}_{m+1}(s)-\tilde{u}_{m}(s)\right) \phi_{1} d W(s)\right) \\
 & -\frac{1}{2} \int_{0}^{t} S(t-s)\left(\left(\tilde{u}_{m+1}(s)-\tilde{u}_{m}(s)\right) F_{\Phi_1}\right) d s \\
= & I_{1}(\omega, t, x)+I_{2}(\omega, t, x)+I_{3}(\omega, t, x)
\end{align*}
For each $\omega \in \Omega$, by Strichartz estimate, similar with the estimate for $A_2$, we have
\begin{equation}\label{I1}
  \left.
   \begin{aligned}
 & \|I_1\|_{\mathcal{X}_{T}}\\
\lesssim &  \| 1_{[0,\tau]} \left( \left(|\lambda_{11}| \left|\tilde{u}_{m+1}(s)\right|^{2\sigma}+|\lambda_{12}|\left|\tilde{v}_{m+1}(s)\right|^{\sigma+1}\left|\tilde{u}_{m+1}(s)\right|^{\sigma-1}\right) \tilde{u}_{m+1}(s) \right.\\
  & \quad \quad \quad  - \left( |\lambda_{11}| \left|\tilde{u}_{m}(s)\right|^{2\sigma}+|\lambda_{12}|\left|\tilde{v}_{m}(s)\right|^{\sigma+1}\left|\tilde{u}_{m}(s)\right|^{\sigma-1}\right) \tilde{u}_{m}(s)  \|_{L^{r^{\prime}}\left([0, t^R_1] ; L^{\frac{2\sigma +2}{2\sigma +1}}\left(\mathbb{R}^N \right)\right)}\\
\lesssim & T^{\frac{2\sigma}{\eta}}( \| \tilde{u}_{m+1}-\tilde{u}_{m} \|_{\mathcal{X}_{T}}+\| \tilde{v}_{m+1}-\tilde{v}_{m} \|_{\mathcal{Y}_{T}} )
   \end{aligned}
   \right.
\end{equation}
We estimate $I_2$ and $I_3$ as same as term $C$ and $D$ in (\ref{ECNLTI}). For $\tilde{v}_{m+1},\tilde{v}_{m}$, we can estimate by same way. Then we finally get
\begin{equation}
  \left.
   \begin{aligned}
      & \| (\tilde{u}_{m+1},\tilde{v}_{m+1})-(\tilde{u}_{m},\tilde{v}_{m})\|_{L^{\rho}(\Omega:(\mathcal{Y}_{T})^2)}\\
    \leqslant & C(R,T,\|\phi_1\|_{W^{1,\infty}},\|\phi_2\|_{W^{1,\infty}})T^{\nu} \| (\tilde{u}_{m+1},\tilde{v}_{m+1})-(\tilde{u}_{m},\tilde{v}_{m})\|_{L^{\rho}(\Omega:(\mathcal{X}_{T})^2)}.
   \end{aligned}
   \right.
\end{equation}
Thus we have $(\tilde{u}_{m+1},\tilde{v}_{m+1})=(\tilde{u}_{m},\tilde{v}_{m})$ on $[0,T]$ for a.e. $\omega \in \Omega$, provided that $T$ is small sufficiently. Hence we can obtain that $(\tilde{u}_{m+1},\tilde{v}_{m+1})=(\tilde{u}_{m},\tilde{v}_{m})$ on $[0,\tau]$ for a.e. $\omega \in \Omega$. $\square$

From Lemma \ref{truncated}, we can define $(u,v)=(u^m,v^m)$ as a local solution to equation (\ref{CNLE}) on $[0,\tau_m]$. Moreover, since $\tau_m$ is increasing with $m$, we can define a stopping time $\tau^{\ast}(u_0,v_0)$ by
\begin{equation*}
\tau^{\ast}(u_0,v_0)=\lim_{m\rightarrow \infty} \tau_m.
\end{equation*}
Then for any stopping time $\tau \leqslant \tau^{\ast}$, equation (\ref{CNLETI}) has a unique $H^1$-solution $(u,v)$ staring from $(u_0,v_0)$ which is almost surely in $(\mathcal{X}_{\tau (\omega)})^2 $.

Now we turn to show that if $\tau^{\ast}(u_0,v_0)<\infty$ then
\begin{equation}\label{H1blowup}
\lim_{t \nearrow \tau^{\ast}(u_0,v_0,\omega)}\sup_{s \leqslant t}(\|u(t))\|+\|v(t))\|)=\infty .
\end{equation}
From the definition of $(u,v)$, if $\tau^{\ast}(u_0,v_0)<\infty$, then we have
\begin{equation}\label{NormInf}
\lim_{t \nearrow \tau^{\ast}(u_0,v_0)} \left( (u,v)  \right)_{(\mathcal{X}_t)^2} = \infty. 
\end{equation}
We prove (\ref{H1blowup}) by contradiction. Let us define
\begin{equation*}
\tilde{\tau}_{R}(\omega)=\inf \left\{t \in\left[0, \tau^{*}\left(u_{0},v_{0}\right)\right), \|u(t)\|_{H^{1}\left(\mathbb{R}^{N}\right)}+\|v(t)\|_{H^{1}\left(\mathbb{R}^{N}\right)} \geq R\right\},
\end{equation*}
then $\tilde{\tau}_{R}$ is a stopping time which satisfies the assumption of Lemma \ref{truncated}. Above estimate (\ref{ET}) claims that
\begin{align*}
 & \mathbb{E} \left( \|u\|_{L^{r}\left(0, \tilde{r}_{R}, W^{1,2\sigma +2}\right)} + \|u\|_{L^{r}\left(0, \tilde{r}_{R}, W^{1,2\sigma +2}\right)} \right) \\ \nonumber
  \lesssim & \|u_0\|_{H^1}+\|v_0\|_{H^1}+ T^{\frac{2\sigma}{\eta}}R^{2\sigma+1} + \left( T^{1-\frac{2}{r}} + T^{\frac{1}{2}} \right) \left( \|u\|_{L^{\rho}(\Omega;\mathcal{X}_T)}+\| v \|_{L^{\rho}(\Omega;\mathcal{X}_T)} \right)
\end{align*}
Assume that
\begin{equation*}
\mathbb{P}\left(\sup _{s \leq \tau^{*}\left(u_{0},v_0 \right)} \left( \|u(s)\|_{H^{1}}+\|v(s)\|_{H^{1}} \right) <+\infty  \text { and }  \tau^{*}\left(u_{0},v_{0}\right)<+\infty\right)>0,
\end{equation*}
then for $R$ sufficiently large, $\mathbb{P}\left( \tilde{\tau}_{R}=\tau^{*}\left(u_{0},v_{0} \right) \right)>0$, which is in contradiction with (\ref{NormInf}). Hence the proof of Theorem \ref{LocalE} is completed. $\square $

\subsection{Global Existence}
In this subsection, we prove the global existence for the mass subcritical case($\sigma < \frac{2}{N} $), the mass critical case($\sigma = \frac{2}{N}$), and the defocusing case.
First of all, we need to estimate the mass and energy of the stochastic Schr\"{o}dinger system (\ref{CNLE}). In order to use It\^o's formula in the following estimation, we need the truncation argument as follows.
For each $k\in \mathbb{N}$, we define the operator $\Theta_k$ by its Fourier transform as
\begin{equation*}
\mathbb{F}(\Theta_k v)(\xi)=\theta_k(|\xi|)\widehat{v}(\xi), \quad t\in\mathbb{R}, \quad \xi\in\mathbb{R}^N,
\end{equation*}
where $\theta_k$ is defined as above, and denote by $(S_k(t))_{t\in\mathbb{R}}$ the linear group by $S_k(t)=\Theta_k S(t)$ or equivalently
\begin{equation*}
\mathbb{F}(S_k(t)v)(\xi)=\theta_k(|\xi|)e^{it|\xi|^2}\widehat{v}(\xi), \quad \xi\in\mathbb{R}^N.
\end{equation*}
Then $S_k(t)$ strongly converges to $S(t)$ in $H^1(\mathbb{R}^N)$.

\begin{lemma}\label{PropM}
Let $(u,v)$ is the solution of (\ref{CNLE}) given by Theorem {\ref{LocalE}} with $u(0)=u_0$, $v(0)=v_0$. Then for any stopping time $\tau< \tau^{\ast}(u_0,v_0)$,
\begin{equation}\label{MassEST}
M(u(\tau))=M(u_0), M(v(\tau))=M(v_0) \quad  \quad a.s.,
\end{equation}
and
\begin{equation}\label{EnergyEST}
\begin{aligned}
 &H\left(u(\tau),v(\tau)\right)\\
=& H\left(u_{0},v_{0}\right)-\operatorname{Im}  \sum_{k \in \mathbb{N}}\int_{\mathbb{R}^{N}} \int_{0}^{\tau} \overline{u} \nabla u \cdot \nabla\left(\phi_1 e_{k}\right) + \overline{v} \nabla v \cdot \nabla\left(\phi_2 e_{k}\right) d x d B_{k}(s)\\
&+\frac{1}{2} \int_{0}^{\tau} \int_{\mathbb{R}^{N}} |u|^2 F_{\phi_1}+|v|^2 F_{\phi_2}+2|uv|\left( F_{\phi_1}F_{\phi_2}\right)^{\frac{1}{2}} d x d s \quad a.s..
\end{aligned}
\end{equation}
\end{lemma}

\noindent{\bf Proof.}
We prove this Lemma by using the following sequence  approximations:
\begin{equation}\label{CNLETIAPP}
  \left\{
   \begin{aligned}
   i d u_{m}^{R}+ & \left(\Theta_{m_{1}} \Delta u_{m}^{R}+\theta_R(\|(u^R_m,v^R_m)\|_{(H^1)^2}) \Theta_{m_{2}}\left[ \left( \lambda_{11}\left|u_{m}^{R}\right|^{2 \sigma}+\lambda_{12}\left|v_{m}^{R}\right|^{\sigma+1}\left|u_{m}^{R}\right|^{\sigma-1} \right) u_{m}^{R} \right]\right) d t\\
        &=u_{m}^{R} \Theta_{m_{2}} \phi_1 d W-\frac{i}{2} u_{m}^{R} F_{\phi_1^{m_{2}}} d t,\\
   i d v_{m}^{R}+ & \left(\Theta_{m_{1}} \Delta v_{m}^{R}+\theta_R(\|(u^R_m,v^R_m)\|_{(H^1)^2}) \Theta_{m_{2}}\left[ \left( \lambda_{21}\left|v_{m}^{R}\right|^{ \sigma-1}\left|u_{m}^{R}\right|^{\sigma+1}+\lambda_{22}\left|u_{m}^{R}\right|^{2 \sigma} \right) v_{m}^{R} \right]\right) d t\\
        &=v_{m}^{R} \Theta_{m_{2}} \phi_2 d W-\frac{i}{2} v_{m}^{R} F_{\phi_2^{m_{2}}} d t.
   \end{aligned}
   \right.
\end{equation}
where index $m=(m_1,m_2)\in\mathbb{N}^2$, $\phi_i^{m_{2}}=\Theta_{m_{2}} \phi_i$ and $F_i^{\phi_{m_{2}}}=\sum_{k=0}^{\infty}\left(\phi_i^{m_{2}} e_{k}(x)\right)^{2}$, $i=1,2$. Since the non-linear term in this system is globally Lipschitz, there exists a unique solution $(u^R_m,v^R_m)$ defined for $t\geqslant 0 $ with initial value $u^R_m(0)=u_0, v^R_m(0)=u_0$. Applying It\^{o}'s formula to $M\left(u_{m}^{R}(t)\right)$, we get
\begin{align*}
  & d \|u^R_M(t)\|^2_{L^2} \\
= & 2 \left( u^R_M (t),i \Theta_{m_{1}} \Delta u^R_M(t) \right)dt \\
  & + 2\left( u^R_M (t), i \theta_R(\|(u^R_m,v^R_m)\|_{(H^1)^2}) \Theta_{m_{2}}\left[ \left( \lambda_{11}\left|u_{m}^{R}\right|^{2 \sigma}+\lambda_{12}\left|v_{m}^{R}\right|^{\sigma+1}\left|u_{m}^{R}\right|^{\sigma-1} \right) u_{m}^{R} \right] \right)dt \\
  & +\left( 2u^R_M (t), -iu^R_M(t)\Theta_{m_{2}} \right)\phi_1^{m_{2}} dW-\frac{1}{2} \left( 2u^R_M (t), u^R_M(t) F_{\phi_i^{m_{2}}} \right)dt \\
  & + \operatorname{Tr} \left( 1, \left( iu^R_M (t) \phi_1^{m_{2}}\right)\left( iu^R_M (t) \phi_1^{m_{2}}\right)^{\ast} \right)dt.
\end{align*}
By the definition of inner product, for any $v \in L^2(\mathbb{R}^N)$, $\left( u^R_M (t),v \right)=0$ if $v$ is a real valued function multiplied by $iu(t,x)$. So the stochastic integration term
\begin{equation*}
  \left( 2u^R_M (t), -iu^R_M(t)\Theta_{m_{2}} \right)\phi_1^{m_{2}}dW=0.
\end{equation*}
Then by integration by part, and taking the imaginary part, we get
\begin{equation*}
  \left( u^R_M (t),i \Theta_{m_{1}} \Delta u^R_M(t) \right) = 0
\end{equation*}
Moreover, since
\begin{equation*}
\left( u^R_M (t), u^R_M(t) F_{\phi_1^{m_{2}}} \right)=\operatorname{Tr} \left( 1, \left( iu^R_M (t) \phi_1^{m_{2}}\right)\left( iu^R_M (t) \phi_1^{m_{2}}\right)^{\ast} \right),
\end{equation*}
it implies that
\begin{align*}
  & d \|u^R_M(t)\|^2_{L^2} \\
= &  2\left( u^R_M (t), i \theta_R(\|(u^R_m,v^R_m)\|_{(H^1)^2}) \Theta_{m_{2}}\left[ \left( \lambda_{11}\left|u_{m}^{R}\right|^{2 \sigma}+\lambda_{12}\left|v_{m}^{R}\right|^{\sigma+1}\left|u_{m}^{R}\right|^{\sigma-1} \right) u_{m}^{R} \right] \right)dt.
\end{align*}
Then by letting $m \rightarrow \infty$ and $R$ sufficient large, we obtain that $d \|u^R_M(t)\|^2_{L^2}=0$ and $\|u(t)\|_{L^2}=\|u_0\|_{L^2}$.
By same way, we obtain $\|v(t)\|_{L^2}=\|v_0\|_{L^2}$. Therefore the system have the conservation of mass (\ref{MassEST}).

Similarly, by applying It\^{o} formula to $H\left(u_{m}^{R}(t)\right)$ and letting $m \rightarrow \infty$, we get
\begin{align*}
&H\left(u^{R}(\tau),v^{R}(\tau)\right) \\
=& H\left(u_{0},v_{0}\right)+\operatorname{Im} \int_{\mathbb{R}^{N}} \int_{0}^{\tau} \overline{u}^{R} \nabla u^{R} \cdot \nabla \phi_1^{m_{2}}d W d x+\operatorname{Im} \int_{\mathbb{R}^{N}} \int_{0}^{\tau} \overline{v}^{R} \nabla v^{R} \cdot \nabla \phi_1^{m_{2}}d W d x \\
&+\operatorname{Im} \int_{0}^{\tau} \int_{\mathbb{R}^{N}} \left(1-\theta_R(\|(u^R,v^R)\|_{(H^1)^2}\right)\left(\left( \lambda_{11}\left|u^{R}\right|^{2 \sigma}+\lambda_{12}\left|v^{R}\right|^{\sigma+1}\left|u^{R}\right|^{\sigma-1} \right) u^{R}\right) \Delta \overline{u}^{R} d x d s \\
&+\operatorname{Im} \int_{0}^{\tau} \int_{\mathbb{R}^{N}} \left(1-\theta_R(\|(u^R,v^R)\|_{(H^1)^2}\right)\left( \left( \lambda_{21}\left|v^{R}\right|^{\sigma-1}\left|u^{R}\right|^{\sigma+1}+\lambda_{22}\left|u^{R}\right|^{2 \sigma} \right) v^{R} \right) \Delta \overline{v}^{R} d x d s \\
&+\frac{1}{2} \int_{0}^{\tau} \int_{\mathbb{R}^{N}}\left|u^{R}\right|^{2}F_{\phi_1} d x d s +\frac{1}{2} \int_{0}^{\tau} \int_{\mathbb{R}^{N}}\left|v^{R}\right|^{2}F_{\phi_2} d x d s+  \int_{0}^{\tau} \int_{\mathbb{R}^{N}}\left|u^{R}v^{R}\right| \left( F_{\phi_1}F_{\phi_2}\right)^{\frac{1}{2}} d x d s.
\end{align*}
By choosing $R$ is sufficiently large, we get the energy estimate (\ref{EnergyEST}), and the proof is completed . $  \square $

\noindent{\bf Proof of Theorem \ref{GlobalE}.}
Let $\tau \leqslant T_0 \wedge \tau^{\ast}(u_0,v_0)$ be a stopping time. For $R>0$, consider the stopping time $\eta_R=inf\left\{t \in[0, \tau^{\ast}(u_0,v_0)],\|(u, v)\|_{\left(\mathcal{X}_{t}\right)^{2}} \leqslant R\right\}$. Then by the conservation of mass (\ref{MassEST}) and the energy estimate (\ref{EnergyEST}), we have
\begin{align*}
& \mathbb{E} \left(\sup _{t \leq \tau \wedge \eta_R}|H(u(t),v(t))|^{2}\right) \\
\leqslant & 2 \mathbb{E}\left(\left|H\left(u_{0},v_0\right)\right|^{2}\right)\\
& +2 \mathbb{E}\left(\sup _{t \leqslant \tau \wedge \eta_R}\left|\int_{\mathbb{R}^{N}} \int_{0}^{t} \overline{u} \nabla u \cdot \nabla \phi_1 d W d x\right|^{2}\right)+2 \mathbb{E}\left(\sup _{t \leqslant \tau \wedge \eta_R}\left|\int_{\mathbb{R}^{N}} \int_{0}^{t} \overline{v} \nabla v \cdot \nabla \phi_2 d W d x\right|^{2}\right)\\
& +\mathbb{E}\left(\int_{0}^{\tau \Lambda \eta_R} \int_{\mathbb{R}^{N}}|u|^{2}F_{\phi_1} d x d s\right)^{2}+\mathbb{E}\left( \int_{0}^{\tau \Lambda \eta_R} \int_{\mathbb{R}^{N}}|v|^{2}F_{\phi_2} d x d s\right)^{2}\\
& +2 \mathbb{E}\left( \int_{0}^{\tau \Lambda \eta_R} \int_{\mathbb{R}^{N}}|uv|\left( F_{\phi_1}F_{\phi_2}\right)^{\frac{1}{2}} d x d s\right)^{2}\\
\leqslant & 2 \mathbb{E}\left(\left|H\left(u_{0},v_0\right)\right|^{2}\right)\\
& +2 \mathbb{E}\left(\int_{0}^{\tau \wedge \eta_R} \int_{\mathbb{R}^{N}} \sum_{k \in \mathbb{N}}\left(\overline{u} \nabla u, \nabla \phi_1 e_{k}\right)^{2} dxd s\right)+2 \mathbb{E}\left(\int_{0}^{\tau \wedge \eta_R}\int_{\mathbb{R}^{N}} \sum_{k \in \mathbb{N}}\left(\overline{v} \nabla v, \nabla \phi_2 e_{k}\right)^{2} dxd s\right) \\
& +2 \mathbb{E}\left( \int_{0}^{\tau \Lambda \eta_R} \int_{\mathbb{R}^{N}}|u|^{2}F_{\phi_1} d x d s\right)^{2}+2\mathbb{E}\left( \int_{0}^{\tau \Lambda \eta_R} \int_{\mathbb{R}^{N}}|v|^{2}F_{\phi_2} d x d s\right)^{2}\\
\leqslant &2 \mathbb{E}\left(\left|H\left(u_{0},v_0\right)\right|^{2}\right)+ C\left(\|F_{\phi_1}\|_{W^{1,\infty}}+\|F_{\phi_2}\|_{W^{1,\infty}} \right) \left( \mathbb{E} \int_{0}^{\tau \Lambda \eta_R}|\nabla u|_{L^{2}}^{4} + |\nabla v|_{L^{2}}^{4} d s \right).
\end{align*}
For the defocusing case, since $\Lambda$ is nonpositive, we have
\begin{align*}
 & \int_{\mathbb{R}^N}|\nabla u(t)|^2+|\nabla v(t)|^2 dx \\
= & 2 H\left( u(t),v(t) \right)+\frac{1}{1+\sigma} \int_{\mathbb{R}^N} \lambda_{11} |u|^{2+2\sigma} +\lambda_{22}|v|^{2+2\sigma}+2\lambda_{12}|v|^{\sigma+1}|u|^{\sigma+1} dx \\
\leqslant & 2 H\left( u(t),v(t) \right).
\end{align*}
Thus by Gronwall inequality and the conservation of mass, we get that
\begin{equation}
\mathbb{E}\left(\sup _{t \leq \tau \wedge \eta_{R}} \left( \|u(t)\|_{H^{1}\left(\mathbb{R}^N\right)}^2 + \|v(t)\|_{H^{1}\left(\mathbb{R}^N \right)}^{2} \right) \right) \leqslant L(\phi_1,\phi_2, u_0, v_0, T_0),
\end{equation}
here $L$ is independent of $R$. Then let $R \rightarrow \infty$, we obtain
\begin{equation}
\mathbb{E}\left(\sup _{t \leq \tau } \left( \|u(t)\|_{H^{1}\left(\mathbb{R}^N\right)}^2 + \|v(t)\|_{H^{1}\left(\mathbb{R}^N \right)}^{2} \right) \right) \leqslant L(\phi_1,\phi_2, u_0, v_0, T_0).
\end{equation}
Thus in the defocusing case, the solution of (\ref{CNLE}) is global.

For the mass critical case $\sigma = \frac{2}{N}$, if $\Lambda$ is positive, then $\lambda_{11}, \lambda_{22}> 0$ and $\lambda_{11} \lambda_{22}-\lambda_{12}^2>0$. By Gagliardo-Nirenberg inequality and the conservation of mass, we have
\begin{align*}
&\int_{\mathbb{R}^N} \lambda_{11} |u|^{2+2\sigma} +\lambda_{22}|v|^{2+2\sigma}+2\lambda_{12}|v|^{\sigma+1}|u|^{\sigma+1} dx \\
=& \|\lambda_{11}^{\frac{1}{4}} |u|^4\|^4_{L^4}+\|\lambda_{22}^{\frac{1}{4}} |v|^4\|^4_{L^4}+\frac{\lambda_{12}}{\sqrt{\lambda_{11}\lambda_{22}}}\|\lambda_{11}^{\frac{1}{4}} uv \|^2_{L^2}\\
\leqslant & K_{opt} \left( \sqrt{\lambda_{11}} \|u\|^2_{L^2}+\sqrt{\lambda_{22}} \|v\|^2_{L^2} \right)\left( \sqrt{\lambda_{11}} \| \nabla u\|^2_{L^2}+\sqrt{\lambda_{22}} \| \nabla v\|^2_{L^2} \right).
\end{align*}
Let $\lambda = \max \lbrace \lambda_{11},\lambda_{22} \rbrace$, then by the conservation of mass,
\begin{equation*}
\frac{1}{2} \left( \|\nabla u(t)\|^2_{L^2}+\|\nabla v(t)\|^2_{L^2} \right)\left[ 1-\frac{\lambda K_{opt}}{2}\left( \sqrt{\lambda_{11}} \|u_0\|^2_{L^2}+\sqrt{\lambda_{22}} \|v_0 \|^2_{L^2} \right) \right] \leqslant H\left( u(t),v(t) \right)
\end{equation*}
Thus when
\begin{equation*}
\left( \sqrt{\lambda_{11}} \|u_0\|^2_{L^2}+\sqrt{\lambda_{22}} \|v_0\|^2_{L^2} \right)<\frac{2}{\lambda K_{opt}},
\end{equation*}
we have
\begin{equation*}
\int_{\mathbb{R}^N}|\nabla u(t)|^2+|\nabla v(t)|^2 dx \lesssim H\left( u(t),v(t) \right).
\end{equation*}
It follows that
\begin{equation*}
\mathbb{E}\left(\sup _{t \leq \tau } \left( \|u(t)\|_{H^{1}\left(\mathbb{R}^N\right)}^2 + \|v(t)\|_{H^{1}\left(\mathbb{R}^N \right)}^{2} \right) \right) \leqslant L(\phi_1,\phi_2, u_0, v_0, T_0).
\end{equation*}

For the mass subcritical case $\sigma < \frac{2}{N}$, by Gagliardo-Nirenberg inequality and the conservation of mass, we get that
\begin{equation*}
\int_{\mathbb{R}^N}|\nabla u(t)|^2+|\nabla v(t)|^2 dx \lesssim H\left( u(t),v(t) \right).
\end{equation*}
It follows that
\begin{equation*}
\mathbb{E}\left(\sup _{t \leq \tau } \left( \|u(t)\|_{H^{1}\left(\mathbb{R}^N\right)}^2 + \|v(t)\|_{H^{1}\left(\mathbb{R}^N \right)}^{2} \right) \right) \leqslant L(\phi_1,\phi_2, u_0, v_0, T_0).
\end{equation*}
Thus in the mass subcritical case, the solution of (\ref{CNLE}) is global. And the proof of Theorem \ref{GlobalE} is completed.

\section{Blow-up solution}

In this section, we prove a general virial identity for the stochastic coupled nonlinear Schr\"{o}dinger system (\ref{CNLE}). As applications, we give a sharp criteria for the blow-up solution of (\ref{CNLE}) when $\frac{2}{N}\leqslant \sigma < \frac{2}{(N-2)^{+}}$.

\begin{lemma}\label{virial}
Let the assumptions of Theorem \ref{LocalE} hold, and $u_0,v_0 \in \Sigma \ a.s.$. Then for any stopping time $ \tau < \tau^{\ast}(u_0,v_0) $, the solution $(u,v) \in C \left( [0,\tau];\Sigma^2 \right) \ a.s. $, the solution satisfies the virial identity
\begin{equation}\label{V}
V(u(\tau),v(\tau))=V(u_0,v_0)+4 \int_0^{\tau} G(u(s),v(s))ds.
\end{equation}
Moreover, for for any stopping time $ \tau < \tau^{\ast}(u_0,v_0) $,
\begin{equation}\label{G}
\begin{aligned}
& G(u(\tau),v(\tau))\\
=& G(u(0),v(0))+4 \int_{0}^{\tau} H(u(s),v(s)) d s \\
&+\frac{2-\sigma N}{\sigma+1} \int_{0}^{\tau}\int_{\mathbb{R}^{N}} \lambda_{11}|u(s)|^{2 \sigma+2}+\lambda_{22}|v(s)|^{2 \sigma+2}+2\lambda_{21}|u(s)|^{ \sigma+1}|v(s)|^{ \sigma+1} dxds \\
&+\sum_{k \in \mathbb{N}} \int_{0}^{\tau} \int_{\mathbb{R}^{N}}|u(s, x)|^{2}x \cdot\nabla\left(\phi_1 e_{k}\right)(x) +|v(s, x)|^{2} x \cdot\nabla\left(\phi_2 e_{k}\right)(x) d x d B_{k}(s).
\end{aligned}
\end{equation}
\end{lemma}

\noindent{\bf Proof.}
We define the truncated form of the variance by
\begin{equation}
V_{\varepsilon}(u,v)=\int_{\mathbb{R}^{N}}e^{-\varepsilon |x|^2}|x|^{2} \left( |u(x)|^{2}+|v(x)|^{2} \right) d x, \quad u,v \in L^2 \left(\mathbb{R}^{N} \right) .
\end{equation}
In order to applying It\^o formula to $V_{\varepsilon}(u,v)$, we need a regularization argument for $u(t),v(t)$. The argument is same as (\ref{CNLETIAPP}) in Lemma \ref{PropM}. Thus we omit it. Then by It\^o's formula, we have
\begin{equation}
\begin{aligned}
 &d V_{\varepsilon}(u(t),v(t)) \\
=&\left( \left( V_{\varepsilon}(u(t),v(t))\right)_u ,i \Delta u(t) \right) d t + \left( \left( V_{\varepsilon}(u(t),v(t))\right)_v ,i \Delta v(t) \right) d t \\
&-\frac{1}{2}\left( \left( V_{\varepsilon}(u(t),v(t))\right)_u, u(t) F_{\phi_1}\right) d t-\frac{1}{2}\left( \left( V_{\varepsilon}(u(t),v(t))\right)_v, v(t) F_{\phi_2}\right) d t \\
&+\frac{1}{2} \operatorname{Tr}\left(\left( V_{\varepsilon}(u(t),v(t))\right)_{uu}(i u(t) \phi_1)(i u(t) \phi_1)^{*}\right) d t \\
&+\frac{1}{2} \operatorname{Tr}\left(\left( V_{\varepsilon}(u(t),v(t))\right)_{vv}(i v(t) \phi_2)(i v(t) \phi_2)^{*}\right) d t \\
&+ \operatorname{Tr}\left(\left( V_{\varepsilon}(u(t),v(t))\right)_{uv}(i u(t) \phi_1)(i v(t) \phi_2)^{*}\right) d t.
\end{aligned}
\end{equation}
After integrating by parts, we have
\begin{equation*}
\left( \left( V_{\varepsilon}(u(t),v(t))\right)_u, i \Delta u(t)\right)=-4 \operatorname{Im} \int_{\mathbb{R}^{N}} e^{-\varepsilon|x|^{2}}\left(1-\varepsilon|x|^{2}\right)(x \cdot \nabla u(t, x)) \overline{u}(t, x) d x.
\end{equation*}
\begin{equation*}
\left( \left( V_{\varepsilon}(u(t),v(t))\right)_v, i \Delta u(t)\right)=-4 \operatorname{Im} \int_{\mathbb{R}^{N}} e^{-\varepsilon|x|^{2}}\left(1-\varepsilon|x|^{2}\right)(x \cdot \nabla v(t, x)) \overline{v}(t, x) d x.
\end{equation*}
Moreover, we have
\begin{equation*}
\frac{1}{2}\left( \left( V_{\varepsilon}(u(t),v(t))\right)_u, u(t) F_{\phi_1}\right)=\int_{\mathbb{R}^{N}} e^{-\varepsilon|x|^{2}}|x|^{2}|u(t, x)|^{2} F_{\phi_1} d x.
\end{equation*}
\begin{equation*}
\frac{1}{2}\left( \left( V_{\varepsilon}(u(t),v(t))\right)_v, v(t) F_{\phi_2}\right)=\int_{\mathbb{R}^{N}} e^{-\varepsilon|x|^{2}}|x|^{2}|v(t, x)|^{2} F_{\phi_2} d x.
\end{equation*}
\begin{equation*}
\frac{1}{2} \operatorname{Tr}\left(\left( V_{\varepsilon}(u(t),v(t))\right)_{uu}(i u(t) \phi_1)(i u(t) \phi_1)^{*}\right)=\int_{\mathbb{R}^{n}} e^{-\varepsilon|x|^{2}}|x|^{2}|u(t, x)|^{2} F_{\phi_1} d x.
\end{equation*}
\begin{equation*}
\frac{1}{2} \operatorname{Tr}\left(\left( V_{\varepsilon}(u(t),v(t))\right)_{vv}(i v(t) \phi_2)(i v(t) \phi_2)^{*}\right)=\int_{\mathbb{R}^{n}} e^{-\varepsilon|x|^{2}}|x|^{2}|v(t, x)|^{2} F_{\phi_2} d x.
\end{equation*}
which implies
\begin{align*}
&-\frac{1}{2}\left( \left( V_{\varepsilon}(u(t),v(t))\right)_u, u(t) F_{\phi_1}\right) d t-\frac{1}{2}\left( \left( V_{\varepsilon}(u(t),v(t))\right)_v, v(t) F_{\phi_2}\right) d t \\
&+\frac{1}{2} \operatorname{Tr}\left(\left( V_{\varepsilon}(u(t),v(t))\right)_{uu}(i u(t) \phi_1)(i u(t) \phi_1)^{*}\right) d t \\
&+\frac{1}{2} \operatorname{Tr}\left(\left( V_{\varepsilon}(u(t),v(t))\right)_{vv}(i v(t) \phi_2)(i v(t) \phi_2)^{*}\right) d t \\
&=0.
\end{align*}
Since $\left( V_{\varepsilon}(u(t),v(t))\right)_{uv}=0$, we have
\begin{equation*}
\operatorname{Tr}\left(\left( V_{\varepsilon}(u(t),v(t))\right)_{uv}(i u(t) \phi_1)(i v(t) \phi_2)^{*}\right)=0
\end{equation*}
From above computations, we get
\begin{align*}
d V_{\varepsilon}(u(t),v(t))=& 4 \operatorname{Im} \int_{\mathbb{R}^{N}} e^{-\varepsilon|x|^{2}}\left(1-\varepsilon|x|^{2}\right) \\
    & \quad \quad \times \left( (x \cdot \nabla u(t, x)) \overline{u}(t, x)+(x \cdot \nabla v(t, x)) \overline{v}(t, x)\right) d x d t,
\end{align*}
For $t<\tau^{\ast}$, integrating in time from $0$ to $t$ yields that
\begin{equation}
\begin{aligned}
V_{\varepsilon}(u(t),v(t))=&V_{\varepsilon}(u_0,v_0) +4 \operatorname{Im} \int_0^t \int_{\mathbb{R}^{N}} e^{-\varepsilon|x|^{2}}\left(1-\varepsilon|x|^{2}\right) \\
    & \quad \quad \times \left( (x \cdot \nabla u(s, x)) \overline{u}(s, x)+(x \cdot \nabla v(s, x)) \overline{v}(s, x)\right) d x d s.
\end{aligned}
\end{equation}
For $k \in \mathbb{N}$, consider stopping time $\tau_k :=\inf \lbrace  t\in[0,T], \|u(t)\|_{H^1 \left(\mathbb{R}^N\right)}+\|v(t)\|_{ H^1 \left(\mathbb{R}^N\right)}\geqslant k \rbrace$. By Cauchy-Schwarz inequality, we have
\begin{equation*}
\begin{aligned}
V_{\varepsilon}\left(u\left(t \wedge \tau_{k}\right), v\left(t \wedge \tau_{k}\right)\right) \leqslant & V_{\varepsilon}\left(u_{0},v_0 \right)+4 k \int_{0}^{t \wedge \tau_{k}} V_{\varepsilon}^{1 / 2}(u(s),v(s)) d s \\
 \leqslant & V_{\varepsilon}\left(u_{0},v_0 \right)+4 k^2 T + \int_{0}^{t \wedge \tau_{k}} V_{\varepsilon}(u(s),v(s)) d s.
\end{aligned}
\end{equation*}
Then by Gronwall inequality, we get
\begin{equation*}
V_{\varepsilon}\left(u\left(t \wedge \tau_{k}\right), v\left(t \wedge \tau_{k}\right)\right) \leqslant \left(  V_{\varepsilon}\left(u_{0},v_0 \right)+4 k^2 T \right) e^T.
\end{equation*}
Thus for almost surely $\omega$, $u(\omega,t),v(\omega,t) \in L^{\infty}\left([0,\tau_k]; \Sigma \right)$ for any $k \in \mathbb{N}$. Letting $k \rightarrow \infty$, we deduce that for almost surely $\omega$, $u(\omega,t),v(\omega,t) \in L^{\infty}\left([0,\tau]; \Sigma \right)$ for any stopping time $\tau^{\ast}(u_0,v_0)$.

It follows from monotone convergence theorem that by letting $\varepsilon \rightarrow 0$,
\begin{equation*}
V(u(\tau),v(\tau))=V(u_0,v_0)+4 \int_0^{\tau} G(u(s),v(s))ds.
\end{equation*}

Applying It\^o's formula to $G(u(t),v(t))$ gives
\begin{align*}
 &d G(u(t),v(t)) \\
=&\left( \left( G(u(t),v(t))\right)_u ,i \left( \Delta u(t)+ (\lambda_{11}|u(t)|^{2\sigma}+\lambda_{12}|v(t)|^{\sigma+1}|u(t)|^{\sigma-1})u(t) \right)\right) d t \\
&+\left( \left( G(u(t),v(t))\right)_v ,i \left( \Delta v(t)+ (\lambda_{21}|v(t)|^{\sigma -1}|u(t)|^{\sigma+1}+\lambda_{22}|u(t)|^{2\sigma})v(t) \right)\right) d t \\
&-\left( \left( G(u(t),v(t))\right)_u, u(t) \right) \phi_1 d W(t)-\left( \left( G(u(t),v(t))\right)_v, v(t) \right) \phi_2 d W(t) \\
&-\frac{1}{2}\left( \left( G(u(t),v(t))\right)_u, u(t) F_{\phi_1}\right) d t-\frac{1}{2}\left( \left( G(u(t),v(t))\right)_v, v(t) F_{\phi_2}\right) d t \\
&+\frac{1}{2} \operatorname{Tr}\left(\left( G(u(t),v(t))\right)_{uu}(i u(t) \phi_1)(i u(t) \phi_1)^{*}\right) d t \\
&+\frac{1}{2} \operatorname{Tr}\left(\left( G(u(t),v(t))\right)_{vv}(i v(t) \phi_2)(i v(t) \phi_2)^{*}\right) d t \\
&+ \operatorname{Tr}\left(\left( G(u(t),v(t))\right)_{uv}(i u(t) \phi_1)(i v(t) \phi_2)^{*}\right) d t.
\end{align*}
From the definition of $ G $ and $ H $, after integration by part, we have
\begin{align}\label{G1}
&\left( \left( G(u(t),v(t))\right)_u, i (\Delta u(t)+ (\lambda_{11}|u(t)|^{2\sigma}+\lambda_{12}|v(t)|^{\sigma+1}|u(t)|^{\sigma-1})u(t) \right) \nonumber \\
&+\left( \left( G(u(t),v(t))\right)_v, i (\Delta v(t)+ (\lambda_{21}|v(t)|^{\sigma -1}|u(t)|^{\sigma+1}+\lambda_{22}|u(t)|^{2\sigma})v(t) \right) d t \nonumber \\
=&\left(2 i x \cdot \nabla u(t)+ iN u(t),i \left( \Delta u(t)+ (\lambda_{11}|u(t)|^{2\sigma}+\lambda_{12}|v(t)|^{\sigma+1}|u(t)|^{\sigma-1})u(t)\right) \right) \nonumber\\
 &+ \left(2 i x \cdot \nabla v(t)+ iN v(t),i \left( \Delta v(t)+ (\lambda_{21}|v(t)|^{\sigma -1}|u(t)|^{\sigma+1}+\lambda_{22}|u(t)|^{2\sigma})v(t) \right) \right) \nonumber \\
=& 4 H(u(t),v(t)) \nonumber \\
 &+\frac{2-\sigma N}{\sigma+1} \int_{\mathbb{R}^N}\lambda_{11}|u(s)|^{2 \sigma+2}+\lambda_{22}|v(s)|^{2 \sigma+2}+2\lambda_{21}|u(s)|^{ \sigma+1}|v(s)|^{ \sigma+1} dx.
\end{align}
A similar calculation gives
\begin{align}\label{G2}
&\left( \left( G(u(t),v(t))\right)_u, u(t) \right)\phi_1 d W(t)+\left( \left( G(u(t),v(t))\right)_v, v(t) \right) \phi_2 d W(t) \nonumber \\
=& \sum_{k \in \mathbb{N}} \left( 2 i x \cdot \nabla u(t, x)+iNu(t, x)), i \overline{u}(t, x)\left(\phi_1 e_{k}\right)(x) \right) d B_k (t) \nonumber \\
 &+\sum_{k \in \mathbb{N}} \left( 2 i x \cdot \nabla v(t, x)+iNv(t, x)), i \overline{v}(t, x)\left(\phi_2 e_{k}\right)(x) \right) d B_k (t)\nonumber \\
=&\sum_{k \in \mathbb{N}} \int_{\mathbb{R}^{N}}|u(t, x)|^{2} x \cdot \nabla\left(\phi_1 e_{k}\right)(x)+|v(t, x)|^{2} x \cdot \nabla\left(\phi_2 e_{k}\right)(x) d x d B_{k}(s).
\end{align}
We also have
\begin{align}\label{G3}
0=&-\frac{1}{2}\left( \left( G(u(t),v(t))\right)_u, u(t) F_{\phi_1}\right) d t-\frac{1}{2}\left( \left( G(u(t),v(t))\right)_v, v(t) F_{\phi_2}\right) d t\nonumber \\
&+\frac{1}{2} \operatorname{Tr}\left(\left( G(u(t),v(t))\right)_{uu}(i u(t) \phi_1)(i u(t) \phi_1)^{*}\right) d t \nonumber \\
&+\frac{1}{2} \operatorname{Tr}\left(\left( G(u(t),v(t))\right)_{vv}(i v(t) \phi_2)(i v(t) \phi_2)^{*}\right) d t.
\end{align}
From (\ref{G1})-(\ref{G3}), we get
\begin{align*}
& dG(u(s),v(s))\\
=& 4  H(u(s),v(s))+\sum_{k \in \mathbb{N}} \int_{\mathbb{R}^{N}}|u(s, x)|^{2}x \cdot\nabla\left(\phi_1 e_{k}\right)(x) +|v(s, x)|^{2} x \cdot\nabla\left(\phi_2 e_{k}\right)(x) d x d B_{k}(s)  \\
&+\frac{2-\sigma N}{\sigma+1} \int_{\mathbb{R}^{N}} \lambda_{11}|u(s)|^{2 \sigma+2}+\lambda_{22}|v(s)|^{2 \sigma+2}+2\lambda_{21}|u(s)|^{ \sigma+1}|v(s)|^{ \sigma+1} dxds.
\end{align*}
Integrating in time from $0$ to $\tau$ yields (\ref{G}). $\square$

\noindent{\bf Proof of Theorem \ref{BLOWUP}.}
Assume that the conclusion of Theorem {\ref{BLOWUP}} is not true, i.e. $\overline{t}<\tau^{\ast}(u_0,v_0)$ a.s.. Then we can take $\tau=\overline{t}$ as a stopping time in Lemma \ref{virial}. Combining with the energy estimate (\ref{EnergyEST}) and (\ref{G}), we have
\begin{align*}
& V(u(\overline{t}),v(\overline{t})) \\
=&V(u_0,v_0)+4G(u_0,v_0)\overline{t}+8H(u_0,v_0)\overline{t}^2\\
&+4\frac{\sigma N-2}{\sigma+1} \int_{0}^{\overline{t}}(\overline{t}-s)\left( \lambda_{11}|u|^{2 \sigma+2}+\lambda_{22}|v|^{2 \sigma+2}+2\lambda_{21}|u|^{ \sigma+1}|v|^{ \sigma+1} \right)dx d s \\
&+4 \int_{0}^{\tau} (\overline{t}-s)^2 \int_{\mathbb{R}^{n}} |u|^2 F_{\phi_1}+|v|^2 F_{\phi_2}+2|uv|\left( F_{\phi_1}F_{\phi_2} \right)^{\frac{1}{2}} d x d s\\
&+4\sum_{k \in \mathbb{N}} \int_{0}^{\overline{t}}(\overline{t}-s) \int_{\mathbb{R}^{N}}|u(s, x)|^{2} x \cdot \nabla\left(\phi_1 e_{k}\right)(x)+|v(s, x)|^{2} x \cdot \nabla\left(\phi_2 e_{k}\right)(x) d x d B_{k}(s)\\
&-8\operatorname{Im} \sum_{k \in \mathbb{N}}\int_{0}^{\tau}(\overline{t}-s)^2 \int_{\mathbb{R}^{n}} \overline{u} \nabla u \cdot \nabla\left(\phi e_{k}\right)(x) + \overline{v} \nabla v \cdot \nabla\left(\phi_2 e_{k}\right)(x) d x d B_{k}(s) .
\end{align*}

For each $t,r\geqslant 0$, let
\begin{align*}
 &V(t,r)\\
=&V(u_0,v_0)+4G(u_0,v_0)t+8H(u_0,v_0)t^2\\
&+4\frac{\sigma N-2}{\sigma+1} \int_{0}^{t}\int_{\mathbb{R}^N} \lambda_{11}|u|^{2 \sigma+2}+\lambda_{22}|v|^{2 \sigma+2}+2\lambda_{21}|u|^{ \sigma+1}|v|^{ \sigma+1} dxd s \\
&+4 \int_{0}^{\tau} (t-s)^2 \int_{\mathbb{R}^{N}} |u|^2 F_{\phi_1}+|v|^2F_{\phi_2}+2|uv|\left( F_{\phi_1}F_{\phi_2} \right)^{\frac{1}{2}} d x d s\\
&+4\sum_{k \in \mathbb{N}} \int_{0}^{r}(\overline{t}-s) \int_{\mathbb{R}^{N}}|u(s, x)|^{2}x \cdot \nabla\left(\phi_1 e_{k}\right)(x)+|v(s, x)|^{2} x \cdot \nabla\left(\phi_2 e_{k}\right)(x) d x d B_{k}(s)\\
&-8\operatorname{Im} \sum_{k \in \mathbb{N}}\int_{0}^{r}(t-s)^2 \int_{\mathbb{R}^{N}}  \overline{u} \nabla u \cdot \nabla\left(\phi_1 e_{k}\right)(x) + \overline{v} \nabla v \cdot \nabla\left(\phi_2 e_{k}\right)(x) d x d B_{k}(s) \\
:=&V(u_0,v_0)+4G(u_0,v_0)t+8H(u_0,v_0)t^2+V_1+4V_2+4V_3-8\operatorname{Im}V_4 .
\end{align*}
For $k \in \mathbb{N}$, consider the stopping time $\tau_k = \inf \lbrace s\in [0,\overline{t}],\|u\|_{H^1}+\|v\|_{H^1} \geqslant k \rbrace $.
For stochastic term $V_3$, by the conservation of mass and It\^o isometry, assuming $t \leqslant \overline{t}$, we get
\begin{align*}
 & \mathbb{E}\left( V_3 \right)^2 \\
=&\mathbb{E}\left( \sum_{k \in \mathbb{N}} \int_{0}^{\tau_{k}}\left| (t-s) \int_{\mathbb{R}^{N}} |u(s, x)|^2 x \cdot \nabla\left(\phi_1 e_{k}\right)(x) +| v(s, x)|^2 x \cdot \nabla\left(\phi_2 e_{k}\right)(x) d x \right|^{2} d s\right)\\
\leqslant& \min_{i=1,2} \left(\|F_{\phi_i}\|_{L^{\infty}}\right) \mathbb{E}\left(\int_{0}^{\tau_{k}}(t-s)^{2} M\left(u_0, v_0\right) V(u(s)) d s\right)\\
\leqslant& \frac{1}{3} \min_{i=1,2} \left(\|F_{\phi_i}\|_{L^{\infty}}\right) \overline{t}^{3} \sup _{s \in[0, \overline{t}]} \mathbb{E}\left(V\left(u\left(s \wedge \tau_{k}\right)\right) M\left(u_0,v_0 \right)\right).
\end{align*}
Thus the stochastic integral $V_3$ is square integrable. It follows that $\mathbb{E}(V_3)=0$. A similar argument as above yields that $\mathbb{E}(V_4)=0$.

For $V_2$, from the conservation of mass we get that
\begin{align*}
V_2 & =\int_{0}^{t}(t-s)^{2} \int_{\mathbb{R}^N} |u(s, x)|^{2}F_{\phi_1}(x)+|v(s, x)|^{2} F_{\phi_2}(x) dxds \\
    & \leqslant \frac{1}{3} t^{3} \min_{i=1,2} \|F_{\phi_i}\|_{L^{\infty}} M\left(u_0,v_0\right) .
\end{align*}

Since $\sigma \geqslant \frac{2}{N}$ and the coefficient matrix $\Lambda$ is negative, we see that $V_1 \leqslant 0$. Now we get that for any $k \in \mathbb{N}$,
\begin{align*}
 \mathbb{E}\left(V\left(t, \tau_{k}\right)\right) & \leqslant \mathbb{E}\left(V\left(u_0,v_0\right)\right)+4 \mathbb{E}\left(G\left(u_0,v_0\right)\right) t \\
 & \quad \quad +8 \mathbb{E}\left(H\left(u_0,v_0\right)\right) t^{2}+\frac{4}{3} t^{3} \min_{i=1,2} \|F_{\phi_i}\|_{L^{\infty}} \mathbb{E}\left(M\left(u_0,v_0\right)\right) .
\end{align*}
We now choose $t=\overline{t}$. From definition, $t_k \rightarrow \tau$ a.s. as $k\rightarrow + \infty$. Then by Fatou's lemma, we have
\begin{align*}
 \mathbb{E}\left(V\left(\overline{t}, \tau_{k}\right)\right) & \leqslant \mathbb{E}\left(V\left(u_0,v_0\right)\right)+4 \mathbb{E}\left(G\left(u_0,v_0\right)\right) t \\
 & \quad \quad +8 \mathbb{E}\left(H\left(u_0,v_0\right)\right) \overline{t}^{2}+\frac{4}{3} \overline{t}^{3} \min_{i=1,2} \|F_{\phi_i}\|_{L^{\infty}} \mathbb{E}\left(M\left(u_0,v_0\right)\right).
\end{align*}
Since $V\left(u_0,v_0\right)$ is nonnegative, the last inequality contradicts with (\ref{Negative}). Thus, Theorem  \ref{BLOWUP} is proved. $\square$

Note that if the expectation over $\omega$ in (\ref{Negative}) is replaced by the expectation over any $\mathcal{F}_0$-measurable subset of $\omega$, the result stated in Theorem \ref{BLOWUP} is also true. Then for any given $\overline{M}>0$ and $\overline{H}>0$, we define
\begin{equation*}
    \mathcal{V}_{\overline{M}, \overline{H}}=\left\{(u,v) \in \Sigma^2, V(u,v)<\overline{M}, G(u,v)<\overline{M}, \|u\|_{L^{2}}^{2}+\|v\|_{L^{2}}^{2}<\overline{M}. H(u,v)<-\overline{H}\right\}
\end{equation*}
We have following corollary
\begin{corollary}
Under the same assumptions as Theorem \ref{BLOWUP}, for any $\overline{M}>0$ and $\overline{t}>0$, there exists a constant $\overline{H}(\overline{t},\overline{M})>0$ such that $\mathbb{P}\left(u_{0} \in \mathcal{V}_{\overline{M}, \overline{H}}\right)>0$ provided with $\mathbb{P}\left(\tau^{*}\left(u_{0}\right) \leq t\right)>0$.
\end{corollary}
\noindent{\bf Proof.}
Let $\Omega_0 =\left\{\omega \in \Omega, u_{0}(\cdot, \omega) \in \mathcal{V}_{\overline{M}, \overline{H}}\right\}$. Then by taking $\overline{H}$ large enough such that
\begin{equation*}
    \overline{M}+4\overline{t}\overline{M}-8\overline{t}^2\overline{H}+\frac{4}{3}\overline{t}^3\min_{i=1,2}\overline{M}<0,
\end{equation*}
and applying Theorem \ref{BLOWUP} with $(u_{0} 1_{\Omega_{0}}, v_0 1_{\Omega_{0}})$, this corollary is proved. $\square$

\end{document}